\documentclass[a4paper]{article}
\usepackage{amsthm,amsmath,amssymb}
\newtheorem{teor}{Theorem}[section]
\newtheorem{defin}[teor]{Definition}
\newtheorem{lemm}[teor]{Lemma}
\newtheorem{osse}[teor]{Remark}
\newtheorem{prop}[teor]{Proposition}
\newtheorem{defi}[teor]{Definition}
\newtheorem{coro}[teor]{Corollary}
\newtheorem{prob}[teor]{Problem}

\newcommand{\bele}{\begin{lemm}\begin{sl}}
\newcommand{\enle}{\end{sl}\end{lemm}}
\newcommand{\bedef}{\begin{defi}\begin{sl}}
\newcommand{\eddef}{\end{sl}\end{defi}}
\newcommand{\bete}{\begin{teor}\begin{sl}}
\newcommand{\ente}{\end{sl}\end{teor}}
\newcommand{\beos}{\begin{osse}\begin{rm}}
\newcommand{\eddos}{\end{rm}\end{osse}}
\newcommand{\bepr}{\begin{prop}\begin{sl}}
\newcommand{\empr}{\end{sl}\end{prop}}
\newcommand{\bepro}{\begin{prob}\begin{rm}}
\newcommand{\empro}{\end{rm}\end{prob}}
\newcommand{\bede}{\begin{defin}\begin{sl}}
\newcommand{\edde}{\end{sl}\end{defin}}
\newcommand{\beco}{\begin{coro}\begin{sl}}
\newcommand{\enco}{\end{sl}\end{coro}}

\newcommand{\quext}{\quad\text}

\newcommand{\RR}{\mathbb{R}}
\newcommand{\EE}{\mathbb{E}}

\newcommand{\beeq}[1]{\begin{equation}\label{#1}}
\newcommand{\eddeq}{\end{equation}}

\newcommand{\beeqa}[1]{\begin{eqnarray}\label{#1}}
\newcommand{\eddeqa}{\end{eqnarray}}

\newcommand{\beal}[1]{\begin{align}\label{#1}}
\newcommand{\eddal}{\end{align}}

\newcommand{\bespl}[1]{\begin{split}\label{#1}}
\newcommand{\edspl}{\end{split}}

\newcommand{\bega}[1]{\begin{gather}\label{#1}}
\newcommand{\edga}{\end{gather}}

\newcommand{\beeqax}{\begin{eqnarray*}}
\newcommand{\eddeqax}{\end{eqnarray*}}

\def\qed{\ifmmode 
  \else \leavevmode\unskip\penalty9999 \hbox{}\nobreak\hfill
  \fi
  \quad\hbox{\hskip.5em\vrule width.4em height.6em depth.05em\hskip.1em}}
\def\endproofsym{\qed}
\renewenvironment{proof}[1][Proof]{\trivlist\item[\hskip\labelsep{\hskip0pt
    {\normalfont\scshape#1.}\hskip .321429\parindent}]\ignorespaces}
{\endproofsym\endtrivlist}
\def\endnobox{\def\endproofsym{}\end{proof}\def\endproofsym{\qed}}

\newcommand{\no}{\nonumber}

\newcommand{\beeqao}{\begin{eqnarray}\no}
\newcommand{\bealo}{\begin{align}\no}
\newcommand{\besplo}{\begin{split}\no}
\newcommand{\begao}{\begin{gather}\no}

\newcommand{\eps}{\varepsilon}

\newcommand{\duav}[1]{\langle{#1}\rangle}

\newcommand{\duavg}[1]{\left\langle{#1}\right\rangle}

\newcommand{\perogni}{\forall\,}
\newcommand{\esiste}{\exists\,}

\newcommand{\io}{\int_\Omega}

\newcommand{\iga}{\int_\Gamma}

\newcommand{\iTT}{\int_0^T}

\newcommand{\iTo}{\iTT\!\io}

\newcommand{\epsi}{\varepsilon}

\newcommand{\dd}{_\delta}

\newcommand{\OO}{_{\Omega}}

\newcommand{\bn}{\boldsymbol{n}}

\newcommand{\dn}{\partial_{\bn}}

\newcommand{\dmax}{_{\delta,\max}}

\newcommand{\lhs}{left hand side}
\newcommand{\rhs}{right hand side}

\DeclareMathOperator{\dive}{div}
\DeclareMathOperator{\deriv}{d}

\DeclareMathOperator{\sign}{sign}

\newcommand{\HUH}{H^1(0,T;H)}

\newcommand{\HUVp}{H^1(0,T;V')}

\newcommand{\LDH}{L^2(0,T;H)}
\newcommand{\LDV}{L^2(0,T;V)}

\newcommand{\LIV}{L^\infty(0,T;V)}

\newcommand{\LDHD}{L^2(0,T;H^2(\Omega))}

\let\TeXchi\chi
\def\chi{{\setbox0 \hbox{\mathsurround0pt
$\TeXchi$}\hbox{\raise\dp0 \copy0 }}}

\newcommand{\aciapo}{\widehat{a}}
\newcommand{\mciapo}{\widehat{m}}
\newcommand{\mjciapo}{\widehat{mj}}

\newcommand{\calE}{{\mathcal E}}

\newcommand{\zzd}{_{0,\delta}}
\newcommand{\zdd}{_{0,\delta}}

\newcommand{\barO}{\overline{\Omega}}

\newcommand{\baru}{\overline{u}}

\newcommand{\dit}{\deriv\!t}
\newcommand{\dis}{\deriv\!s}

\newcommand{\ddt}{\frac{\deriv\!{}}{\dit}}

\newcommand{\calL}{{\cal L}}

\numberwithin{equation}{section}
\begin{document}

\title{A Cahn-Hilliard equation
  with singular diffusion}

\author{Giulio Schimperna\\
Dipartimento di Matematica, Universit\`a di Pavia,\\
Via Ferrata~1, I-27100 Pavia, Italy\\
E-mail: {\tt giusch04@unipv.it}\\
\and
Irena Paw\l ow\\
Systems Research Institute,\\
Polish Academy of Sciences\\
\mbox{}~~and Institute of Mathematics and Cryptology,\\
Cybernetics Faculty,\\
Military University of Technology,\\
S.~Kaliskiego 2, 00-908 Warsaw, Poland\\ 
E-mail: {\tt Irena.Pawlow@ibspan.waw.pl}
}


\maketitle
\begin{abstract}
 In the present work, we address a class of Cahn-Hilliard equations
 characterized by a singular diffusion term.
 The problem is a simpified version with constant mobility 
 of the Cahn-Hilliard-de Gennes model of phase separation in binary, 
 incompressible, isothermal mixtures of polymer molecules.
 It is proved that for any final time $T$, the problem admits a unique 
 energy type weak solution, defined over $(0,T)$. For any $\tau > 0$
 such solution is classical in the sense of belonging to a suitable
 H\"{o}lder class over $(\tau, T)$, and enjoys the property of 
 being separated from the singular values corresponding to pure phases. 
\end{abstract}

\noindent {\bf Key words:}~~Cahn-Hilliard equation,
singular diffusion, variational formulation, existence theorem.

\vspace{2mm}

\noindent {\bf AMS (MOS) subject clas\-si\-fi\-ca\-tion:}%
~~35K35, 35K67, 35A01, 82D60.

\vspace{2mm}


\section{Introduction}
\label{sec:intro}

This paper is devoted to the mathematical analysis of the following
class of parabolic systems:
\begin{align}\label{CH1}
  & u_t - \Delta w = 0,\\
 \label{CH2}
  & w = - a(u) \Delta u - \frac{a'(u)}2 |\nabla u|^2 
   + f(u) - \lambda u + \epsi u_t,
\end{align}
on $(0,T) \times \Omega$, $\Omega$ being a bounded 
smooth subset of $\RR^3$ and $T>0$ an assigned final
time. The system is coupled with
the initial and boundary conditions
\begin{align}\label{iniz-intro}
  & u|_{t=0} = u_0,
  \quext{in }\,\Omega,\\
 \label{neum-intro}
  & \dn u = \dn w = 0,
  \quext{on }\,\partial\Omega,\
  \quext{for }\,t\in(0,T)
\end{align}
and represents a variant of the Cahn-Hilliard model 
for phase separation in binary materials. 
Here, $\lambda$ and $\epsi$ are nonnegative parameters,
where the case $\epsi>0$ accounts for a viscosity effect
that may appear the frame of Cahn-Hilliard models
(see, e.g., \cite{Gu}). Moreover,
the function $f$ stands for the derivative 
of the so-called logarithmic potential  
\begin{equation}\label{logpot}
  F(r) = (1-r)\log(1-r)+(1+r)\log(1+r),
   \quad r\in[-1,1],
\end{equation}
namely,
\begin{equation}\label{f}
  f(r) = \log(1+r)-\log(1-r)= \log \frac{1+r}{1-r},
   \quad r\in(-1,1).
\end{equation}
Finally, we assume that the function $a$
has the form
\begin{equation}\label{defia}
   a(r) = \frac2{1-r^2} = f'(r).
\end{equation}
The initial-boundary value problem given by \eqref{CH1}-\eqref{iniz-intro}
will be noted as Problem~(P) in the sequel.



The above problem is a simplified version (assuming constant mobility and
normalized physical quantities) of the Cahn-Hilliard-de Gennes model 
\cite{dG80} of phase separation in binary, incompressible, isothermal
mixtures of polymer molecules of different types $i = 1$ and $i = 2$,
quenched below a critical temperature.
 
Following the Cahn-Hilliard-de Gennes theory
\cite{dG80}, \cite{dG85} (see also \cite{Bin83}, \cite{NL84}, \cite{Pin81}),
we describe now briefly
the physical basis of the presented model. 
Each molecule $i$ in the mixture consists of $N_i$ segments of size $\sigma_i$
(so-called lattice constant) with the quantity 
$R^2_{g_i} = \frac{1}{6} N_i \sigma^2_i $ 
denoting the mean radius of gyration of the i-th polymer molecule.
The variables $u$ and $w$ in Problem (P) have the meaning of the rescaled 
order parameter and the exchange chemical potential. For clarity of physical 
description let us use just now in place of $u \in [-1,1]$ the original variable 
$\chi = \chi_1 \in [0,1]$, denoting the volume fraction of component $i = 1$.
Then, by the incompressibility condition, $\chi_1 + \chi_2 = 1$
everywhere in the sample. The variable
\begin{equation}\label{defiu} 
    u:=2\chi - 1 \in [-1,1]
\end{equation}
is introduced for the sake of mathematical convenience.

The model is governed by the Flory-Huggins-de Gennes free energy 
functional \cite{dG80} 
which in the isothermal case is
\begin{equation}\label{defiFHdG} 
\mathcal{F}_{FHdG}(\chi) =\io \Big( F_{FH}(\chi) + \frac{1}{2} a(\chi) |\nabla \chi|^2 \Big).
\end{equation}
The homogeneous (volumetric) free energy $F_{FH}(\chi)$ has the Flory-Huggins form
\begin{equation}\label{defiFH}
  F_{FH}(\chi) = \frac{1}{N_1}\chi \log \chi + \frac{1}{N_2}(1 - \chi) \log(1 - \chi)
      + \lambda \chi(1 - \chi),
\end{equation}
where $\lambda$, called the Flory-Huggins interaction parameter, measures the strength
of interaction between two kinds of species. The second term in the integrand
\eqref{defiFHdG} gives the weighted contribution due to composition gradients.
By de Gennes theory \cite{dG80}, it has the characteristic singular 
form
\begin{equation}\label{apolymer}
  \frac{1}{2}a(\chi) = \frac{1}{6}\Big( \frac{R^2_{g_1}}{N_1 \chi} + \frac{R^2_{g_2}}{N_2 (1 - \chi)} \Big) 
  = \frac{1}{36}\Big( \frac{\sigma_1^2}{\chi} + \frac{\sigma_2^2}{1 - \chi} \Big)      
\end{equation}
which is due to the connectivity of the chains that constitute the polymer molecules.
In the symmetrical case $N_1 = N_2 = N, \; \sigma_1 = \sigma_2 =\sigma,
\; R^2_{g_1} = R^2_{g_2} = R^2_{g} $, the expression \eqref{apolymer} simplifies to
\begin{equation}\label{apolymersymm}
 \frac{1}{2} a(\chi) = \frac{\sigma^2}{36} \frac{1}{\chi (1 - \chi)}.       
\end{equation}
Here we note that the singular form of $a(\chi)$  introduces an infinite 
energy penalty near the pure phases. 
This fact turns out to have an important mathematical consequence
related to the separation property of the solution
(see \eqref{separs} and
Remark~\eqref{mirzel} below). 

Under the incompressibility condition, 
$\chi_1 + \chi_2 = 1$, the conservation of mass for component $1$ is
given by
\begin{equation}\label{massbalance}
  \chi_t + \nabla \cdot j = 0, \quad j = - \Lambda(\chi) \nabla w.  
\end{equation}
Here $j$ is the mass flux defined as the product of the Onsager diffusion coefficient 
$\Lambda(\chi)$ (effective mobility) and the gradient of the exchange potential $w$
which is the difference between the chemical potentials of the components,
$w = w_1 - w_2$. 
Shortly, we shall refer to $\Lambda$ and $w$ as the mobility 
and the chemical potential, respectively.

According to the derivations in \cite{dG80}, \cite{Bin83}, \cite{BJL83}, 
the expression for $\Lambda$
takes the form
\begin{equation}\label{Onsager}
  \Lambda(\chi) = \frac{\Lambda_1 \Lambda_2}{\Lambda_1  +  \Lambda_2}  = \chi (1 - \chi) \Lambda_0,       
\end{equation}
where $\Lambda_i = \chi_i \Lambda_0$ is the Onsager coefficients for the $i$-th component
and $\Lambda_0$ is a positive constant.
 
We point out that the concentration dependence of the mobility $\Lambda$ 
is also typical for phase separation of small molecule systems, described by 
the classical Cahn-Hilliard equation \cite{Ca}, \cite{CH}; 
see also \cite{EG}, \cite{GNC00}, \cite {NB89} and the references therein.

The chemical potential $w$ is defined as the first variation of 
the functional \eqref{defiFHdG} yielding the following equivalent expressions:
\begin{align}\no
  w & = \frac{\delta \mathcal{F}_{FHdG}}{\delta \chi} \\ 
  \nonumber
      & =   F^{\prime}_{FH}(\chi) + \frac{1}{2} a^{\prime}(\chi)|\nabla \chi|^2 - 
  \dive (a(\chi) \nabla \chi) \\
  \label{chempot}
      & =   F^{\prime}_{FH}(\chi) - \frac{1}{2} a^{\prime}(\chi)|\nabla \chi|^2
  - a(\chi) \Delta \chi \\
      & =   F^{\prime}_{FH}(\chi) - \sqrt{a(\chi)} \dive (\sqrt{a(\chi)} \nabla \chi).
      \nonumber
\end{align}
As in the standard Cahn-Hilliard theory, the chemical potential can 
include an additional term $\epsi \chi_t$, with a positive constant $\epsi$, 
accounting for possible viscous effects.

The equations \eqref{massbalance}, \eqref{Onsager}, and \eqref{chempot}
with the term $\epsi \chi_t, \;\eps \geq 0$, lead to the following
degenerate singular Cahn-Hilliard-de Gennes polymer system:
\begin{align}\label{CH1p}
  & \chi_t - \Lambda_0 \dive (\chi (1 - \chi) \nabla w)= 0,\\
 \label{CH2p}
  & w = - a(\chi) \Delta \chi - \frac{a'(\chi)}2 |\nabla \chi|^2 
   + F^{\prime}_{FH}(\chi)+ \epsi \chi_t
\end{align}
to be considered together with the initial 
and appropriate boundart conditions.

As a first step in the study of such system we consider the case 
of constant mobility by linearizing the principal part
$\chi (1-\chi)$ in \eqref{CH1p}. 
Moreover, for mathematical convenience,
we replace $\chi$ by the variable $u$ defined by \eqref{defiu},
restrict ourselves to the symmetrical mixture \eqref{apolymersymm},
and set all physical constants equal to unity.
This leads to  \eqref{CH1}-\eqref{CH2}.

Apart from the well-known Cahn-Hilliard-de Gennes model described above, 
we mention also recently developed two-fluids models for viscoelastic 
phase separation in polymer solutions, see \cite{ZZE06}
and the references therein. 

It is known that the process of polymer phase separation
is important both for its theoretical aspect and due to
unusual morphology for specific materials applications.
On the contrary to extensive physical and numerical literature 
on polymer mixtures (for review see e.g., \cite{Bin83}, \cite{dG85}, 
\cite{NL84}, \cite{Nos87}, \cite{AP96}),
mathematical aspects of the corresponding models have not been 
so far addressed sufficiently.
The structure of steady-state solutions to degenerate singular
polymer system \eqref{CH1p}-\eqref{CH2p} has been analysed by 
Mitlin et al. \cite{MiMaE85}, \cite{MiMa90},
and Witelski \cite{Wit98}. The existence of weak solutions to the 
Cahn-Hilliard equation with logarithmic potential $F(u)$, degenerate mobility
$\Lambda(u) =1-u^2$, and constant gradient coefficient $a>0$ has been studied 
by Elliott and Garcke \cite{EG}; we refer also to \cite{DpGG} and \cite{GNC00} 
for further study of degenerate problems.
For the standard Cahn-Hilliard equation with logarithmic potential $F(u)$
and constant coefficient $a>0$, optimal regularity of weak solutions is 
analyzed in \cite{MZ}, where in particular the separation property~\eqref{separs}
is obtained in space dimensions 1 and 2.
The existence and uniqueness of weak solutions to the Cahn-Hilliard equation 
with logarithmic potential and nonlinear positive, bounded coefficient
$a(u)$ has been recently proved by the authors in \cite{SP11}.

We mention also a closely related sixth order Cahn-Hilliard type problem with
nonlinear coefficient $a(u)$, considered recently in \cite{SP11} 
for a singular (e.g., logarithmic) potential, and in \cite{PZ11}
for a polynomial potential. As a special case, in \cite{SP11} 
the behavior of the solutions when the sixth order term is let tend to zero 
was analysed.

Finally, we point out that the theoretical investigation 
of polymer models was initiated by Alt and the second author 
in \cite{AP96}, where general nonisothermal phase transition 
models with a conserved order parameter
have been derived and, in particular,  polymer free energy models 
have been presented along with an extensive list of references.
The authors of \cite{AP96} have obtained also
some partial, unpublished results  \cite{APNotes96} on the existence 
of weak solutions to degenerate singular polymer model 
\eqref{CH1p}-\eqref{CH2p}  
by applying the methods due to Elliott and Garcke \cite{EG},
and Elliott and Luckhaus \cite{EL91}. 
This unsolved problem has become the motivation of the present study
which uses a different approach developed previously in \cite{SP11}.

As already mentioned, as a first step of the analysis 
we assume that the mobility is constant, 
and extend our methods applied in \cite{SP11}
in the case of a nonlinear (but bounded)
coefficient $a(u)$ to singular $a(u)$.

There are two main ideas behind our approach. 
The first one, standard in the analysis of Cahn-Hilliard systems,
exploits the characteristic variational structure 
of the system \eqref{CH1}-\eqref{CH2}.
The second non-standard one consists in introducing appropriate 
changes of variables. 

The variational structure becomes evident by (formally) 
testing \eqref{CH1} by $w$, \eqref{CH2} by $u_t$,
taking the difference of the obtained relations, 
integrating with respect to space variables,
using the {\sl no-flux}\/ conditions \eqref{neum-intro},
and performing suitable integrations
by parts. Then one readily gets the {\sl a-priori}\/ bound
\begin{equation}\label{energyineq}
  \ddt\calE(u) + \| \nabla w \|_{L^2(\Omega)}^2 
   + \epsi \| u_t \|_{L^2(\Omega)}^2 
   = 0,
\end{equation}
which has the form of an {\sl energy equality}\/ for the 
{\sl energy functional}
\begin{equation}\label{defiE}
  \calE(u)=\io \Big( \frac{a(u)}2 |\nabla u|^2 
   + F(u) - \lambda \frac{u^2}2 \Big),
\end{equation}
where the interface (gradient) part contains the nonlinear 
function $a$. In other words, the system \eqref{CH1}-\eqref{CH2}
arises as the $(H^1)'$-gradient flow problem for the 
functional $\calE$.

However, the energy estimate~\eqref{energyineq} is not sufficient
to obtain existence of a solution to \eqref{CH1}-\eqref{CH2}
via approximation-compactness methods. Actually, to apply 
this strategy, one also 
needs some control of the second space derivatives of $u$ and of the 
singular coefficients $a$ and $a'$ in \eqref{CH2}. To obtain this,
two changes of variables will play an important role. The first one, 
motivated by the structure of the fourth formula in 
\eqref{chempot} and applied also in \cite{SP11}, consists in introducing 
the variable $z$ such that $\nabla z = \frac{1}{\sqrt{2}} \sqrt{a(u)} \nabla u$.
Then, the Laplacean of $z$ can be (formally)
estimated simply by testing the equivalent
$z$-formulation of \eqref{CH2} (namely, \eqref{CH2z} below) by 
$-\Delta z$. This gives the desired control on second space derivatives.

On the other hand, even in the equivalent formulation \eqref{CH2z},
one has to control a coefficient (namely, $\phi'(u)$), that explodes
polynomially fast as $|u|$ approaches $1$. This is a nontrivial
issue since the nonlinear term $f(u)$ in \eqref{CH2} (or in 
\eqref{CH2z}), for which is relatively simple to get
a $L^2$-control, explodes only {\sl logarithmically}\/ fast 
(and, hence, an $L^p$-estimate of $f(u)$ would help to control
$\phi'(u)$, or $a(u)$, only for $p=\infty$). 
To overcome this further difficulty, a second additional
change of variable comes into help. Namely, we set
$v := f(u)$, which represents 
the monotone part of the volumetric 
chemical potential. We then see that the formulation 
of equation \eqref{CH2} in terms of $v$ 
(namely, \eqref{CH2v} below) does no longer contain
singular coefficients of polynomial type;
it presents, however, the cubic term $u |\nabla v|^2$.
To control it, we use techniques based on entropy-type 
estimates (cf.~\cite{DpGG}), with a rather careful and 
{\sl ad-hoc}\/ choice of test functions
(cf.~\eqref{defim} and \eqref{defimp} below).
In this way we can both prove an $L^\infty$-bound 
for $v$ (and consequently the ``separation property''
\eqref{separs}), for strictly positive times,
and also control the cubic term $u |\nabla v|^2$
starting from the initial time $t=0$. This is the key
step that permits to get existence of a weak solution
(for the $v$-formulation of the system) for initial
data $u_0$ that have only the natural energy regularity
(i.e.,~such that $\calE(u_0)<+\infty$).

Regarding additional properties of solutions, we remark
that, as already noticed in \cite{DNS}, the energy
\eqref{defiE} with the coefficient $a$ given by
\eqref{defia} is convex with respect to~$u$
(up to the lower order $\lambda$-perturbation). This 
basic property permits to prove parabolic time-regularization 
properties of weak solutions, as well as uniqueness,
in a relatively standard way.

The plan of the paper is as follows. 
In the next Section~\ref{sec:not} we will present
the main assumptions and give the statement of 
our main results. In Section~\ref{sec:strong}, we will 
prove local existence and uniqueness of a strong (i.e.,
lying in a suitable H\"older class) solution. 
The main a-priori estimates needed in the proof
of global existence will be detailed in the subsequent
Section~\ref{sec:apriori}. On the basis of these
estimates, in the final Section~\ref{sec:weak}
we shall show global existence, uniqueness, and 
time-regularization properties of weak solutions.


\medskip
\noindent
{\bf Acknowledgment.}~~The second author has greatly benefited 
from the cooperation on phase transition models with Prof.~H.W.~Alt
during her stay at the Institute of Applied Mathematics,
University of Bonn. The support of SFB~256
``Nichtlineare partielle Differentialgleichungen'' 
is greatly acknowledged.



\section{Notation and main results}
\label{sec:not}

Let $\Omega$ be a smooth bounded domain of $\RR^3$
of boundary $\Gamma$, $T>0$ a given final time, and
let $Q:=(0,T)\times\Omega$. Let $H:=L^2(\Omega)$,
endowed with the standard scalar product $(\cdot,\cdot)$
and the norm $\| \cdot \|$. Let also $V:=H^1(\Omega)$.
We identify $H$ with $H'$ so that the chain of continuous
embeddings $V\subset H \subset V'$ holds.
We indicate by $\duav{\cdot,\cdot}$ the duality between
$V'$ and $V$ and by $\|\cdot\|_X$ the norm in 
the generic Banach space $X$. We note as $A$ the
weak Laplace operator with no-flux boundary conditions, 
namely
\begin{equation} \label{defiA}
   A: V \to V', \qquad
    \duav{Av,z} := \io \nabla v \cdot \nabla z,
      \quad \perogni v, z \in V.
\end{equation}
We also set
\begin{equation} \label{defiW}
  W := \big\{ z \in V:~Az \in H \big\}
   = \big\{ z \in H^2(\Omega):~\dn z = 0~\text{on }\Gamma \big\},
\end{equation}
which is a closed subspace of $H^2(\Omega)$.
In all what follows we shall assume that $F$, $f$ and 
$a$ are given, respectively, by \eqref{logpot}, 
\eqref{f} and \eqref{defia}. Moreover, we will
assume that $u_0$ is an initial datum having
{\sl finite energy}\, $\calE$ (cf.~\eqref{defiE}),
namely
\begin{equation} \label{EE0}
  \EE_0:= \calE(u_0) < +\infty.
\end{equation}
It is worth noting that, since $F$ has
the expression \eqref{logpot},
the above is equivalent to asking 
\begin{equation} \label{regou0}
  -1 \le u_0 \le 1~~\text{a.e.~in }\,\Omega, \qquad
  \qquad  a^{1/2}(u_0) \nabla u_0  \in H.
\end{equation}
Moreover, due to \eqref{defia}, if $u_0$ satisfies \eqref{regou0}
then $u_0\in V$, in particular. Letting, for a generic
summable function $v:\Omega\to \RR$, $v\OO$ denote its
spatial mean value, we will also assume that
\begin{equation} \label{regou02}
  m := (u_0)\OO = \frac1{|\Omega|} \io u_0
   \in (-1,1).
\end{equation}
In other words, we cannot admit the case when
$u_0$ coincides with $+1$ (or with $-1$) almost everywhere
in~$\Omega$. This is a standard assumption when dealing
with Cahn-Hilliard systems containing constraint terms
(cf, e.g., \cite{KNP} for more details).

Let us note also that, substituting the expression
\eqref{defia} for $a$, \eqref{CH2} can be rewritten as
\begin{equation}\label{CH2a}
  w = - \frac2{1-u^2} \Delta u 
   - \frac{2u}{(1-u^2)^2} |\nabla u|^2
   + f(u) - \lambda u + \epsi u_t.
\end{equation}
It is now convenient to introduce
a couple of additional variables which permit to give
alternative formulations of \eqref{CH2a}.
To start with, we compute some derivatives of $a$. 
From \eqref{logpot} and \eqref{defia}, we have
\begin{equation}\label{aprimo}
   a'(r) = f''(r) = \frac{4r}{(1-r^2)^2},
\end{equation}
as well as 
\begin{equation}\label{asecondo}
   a''(r) = f'''(r) = \frac{4(1+3r^2)}{(1-r^2)^3}.
\end{equation}
In the sequel, for a generic locally integrable 
real-valued function $\psi$ defined in an open 
neighbourhood of $0$, we will write
\begin{equation}\label{deficiapo}
  \widehat{\psi}(r):=\int_0^r \psi(s)\,\dis.
\end{equation}
Then, of course, we have
\begin{equation}\label{aciapo}
   \aciapo(r) = \int_0^r a(s)\,\dis 
    = f(r) = \log (1+r) - \log (1-r).
\end{equation}
Next, we introduce
\begin{equation}\label{defiphi}
  \phi(r) := \frac{\sqrt2}2\int_0^r a^{1/2}(s)\,\dis 
    = \int_0^r \frac{ \dis }{ (1-s^2)^{1/2} }
    = \arcsin r.
\end{equation}
Then, we notice that, setting 
\begin{equation} \label{defiz}
  z := \phi(u) = \arcsin u,
\end{equation}
equation \eqref{CH2} can be rewritten in the equivalent
form 
\begin{equation}\label{CH2z}
  w = - 2 \phi'(u) \Delta z + f(u) 
   - \lambda u + \epsi u_t.
\end{equation}
Next, we put $v:= f(u)$. Then, a simple computation gives
\begin{equation}\label{co61}
  u = f^{-1}(v) = \frac{e^v - 1}{e^v + 1} =: j(v).
\end{equation}
Moreover, \eqref{CH2} can be rewritten as
\begin{equation}\label{CH2v0}
  w = - \Delta v + v + \frac{a'(u)}2 | \nabla u |^2 
   - \lambda u + \epsi u_t.
\end{equation}
Noting that
\begin{equation}\label{co62}
  \frac{a'(u)}2 | \nabla u |^2  
   = \frac{2u}{(1-u^2)^2} | \nabla u |^2
   = \frac{u}2 \big| f'(u)\nabla u \big|^2,
\end{equation}
we finally obtain from \eqref{co61}
\begin{equation}\label{CH2v}
  w = - \Delta v + v + \frac{j(v)}2 | \nabla v |^2  - \lambda j(v) + \epsi u_t
   = - \Delta v + v + \frac{1}2 \frac{e^v - 1}{e^v + 1} | \nabla v |^2 
      - \lambda j(v) + \epsi u_t.
\end{equation}

In the sequel, we shall indicate by $c$ a generic positive
constant, whose value may vary on occurrence,
allowed to depend on the parameters of the system
(more precisely, on the functions $a$ and $f$, on $\lambda$
and on $\Omega$), and, in particular, not on approximating
parameters. Moreover, the constants $c$ will not be allowed
to depend on the choice of initial data. However,
they may depend on the prescribed mean value $m$.
The notation $\kappa$ will be used for positive constants
(depending on the same quantities as $c$) appearing in 
estimates from below. We will also use the notation 
$Q(\cdot)$ (or $Q(\cdot,\cdot)$), 
with $Q$ indicating a computable function with values in
$[0,+\infty)$, increasingly monotone in each of its argument,
whose expression can depend on the same quantities as $c$.
For instance, the expression $Q(\EE_0)$ will stand for a 
monotone function of the ``initial energy'' $\EE_0$.

We can now introduce the concepts of ``classical'' and of ``weak'' 
solution needed in the subsequent analysis. 
%
%
%
\bede\label{def:classsol}
 A ``strong'', or ``classical'', solution to~{\rm Problem~(P)}
 over the time interval $(0,T)$
 is a couple $(u,w)$ with the regularity 
 \begin{align}\label{regoust}
   & u\in W^{1,\infty}(0,T;V') \cap H^1(0,T;V) 
      \cap L^\infty(0,T;H^2(\Omega)),
    \qquad \epsi u\in W^{1,\infty}(0,T;H),\\
  \label{regowst}
   & w\in L^\infty(0,T;V) \cap L^2(0,T;H^3(\Omega)), \qquad 
    \epsi w\in L^{\infty}(0,T;H^2(\Omega)),
 \end{align}
 satisfying, a.e.~in~$(0,T)$, the equations
 \begin{align}\label{CH1s}
   & u_t - \Delta w = 0, \quext{a.e.~in }\,\Omega,\\
  \label{CH2s}
   & w = - a(u) \Delta u - \frac{a'(u)}2 |\nabla u|^2 
   + f(u) - \lambda u + \epsi u_t,
    \quext{a.e.~in }\,\Omega,\\
  \label{neums}
   & \dn u = \dn w = 0, \quext{a.e.~on }\,\Gamma,
\end{align}
together with the initial condition $u|_{t=0} = u_0$
and, for all $(t,x)\in [0,T]\times \barO$,
the\/ {\rm separation property}
\begin{equation} \label{separs}
  -1 + \epsilon \le u(t,x) \le 1 - \epsilon,
   \quext{for some }\,\epsilon>0.
\end{equation}
\edde
\beos\label{morereg}
 Thanks to \eqref{separs}, the component $u$ of 
 any ``classical'' solution to Problem~(P) is uniformly
 separated from the singular values $\pm1$ of $a$ and $f$. 
 Hence, by applying the standard theory of quasilinear 
 parabolic equation and a bootstrap argument, 
 we can see that $(u,w)$ is in fact smoother, with 
 its regularity being limited only by the regularity of the
 initial datum. In other words, at 
 least for times $t>0$, a classical solution
 can be thought to be arbitrarily regular.
\eddos
\bede\label{def:weakssol}
 A ``weak'', or ``energy'', solution to~{\rm Problem~(P)}
 over the time interval $(0,T)$
 is a couple $(u,w)$ with the regularity 
 \begin{align}\label{regou4}
   & u\in \HUVp\cap L^\infty(0,T;V) \cap L^\infty((0,T)\times\Omega),
    \qquad \epsi u\in \HUH,\\
   \label{regoFu4}
   & F(u) \in L^\infty(0,T;L^1(\Omega)),\\
   \label{regofu4}
   & v = f(u) \in L^2(0,T;V),
    \qquad \nabla v\in L^p((0,T)\times\Omega)~~
     \text{for some }\,p>2,\\
  \label{regow4}
   & w\in L^2(0,T;V),
 \end{align}
 satisfying, a.e.~in~$(0,T)$, the equations
 \begin{align}\label{CH1w}
   & u_t + A w = 0, \quext{in }\,V',\\
  \label{CH2w}
   & w = - \Delta v + v + \frac{u}2 | \nabla v |^2  - \lambda u + \epsi u_t,
    \quext{a.e.~in }\,\Omega,\\
  \label{uv}
   & v = f(u), \quext{a.e.~in }\,\Omega,\\
  \label{neumw}
   & \dn v = 0, \quext{a.e.~on }\,\Gamma.
 \end{align}
 together with the initial condition~$u|_{t=0} = u_0$.
\edde
\noindent%

We can now state our main result, regarding existence,
uniqueness, and regularization properties of 
weak solutions to~Problem~(P):
\bete\label{teo:main}
 Let $f$ and $a$ be given by \eqref{f}, \eqref{defia}.
 Let $u_0$ satisfy\/ \eqref{EE0} and~\eqref{regou02}.
 Finally, let $\Omega$ be\/ {\rm convex}.
 Then, for any $T>0$, {\rm Problem~(P)} admits at least
 a weak solution  $(u,w)$ defined over $(0,T)$. 
 Moreover, for any $\tau>0$, $(u, w)$ is a ``classical'' solution
 over $(\tau,T)$. In particular, the separation 
 property~\eqref{separs} holds on $(\tau,T)$
 with
 \begin{equation}\label{separepsi}
   \epsilon^{-1} = Q(\EE_0,\tau^{-1}).
 \end{equation}
 Finally, uniqueness holds in the class of weak 
 solutions that are classical for strictly positive times. 
\ente
\noindent%
For strictly positive times, \eqref{CH2} can be interpreted 
in any of the equivalent formulations \eqref{CH2a}, \eqref{CH2z}, 
or \eqref{CH2v}; actually, $u$ is a classical solution
for $t>0$. On the other hand, when looking at the behavior
near $t=0$, it is crucial to view \eqref{CH2} in 
the form~\eqref{CH2w} which appears to be the only 
formulation permitting to take the approximation limit 
$\delta\searrow 0$ starting from the initial time.
\beos\label{morereg2}
 If we have in addition that the initial datum
 satisfies $u_0\in H^2(\Omega)$ with
 $u_0(x)\in [-1+\epsilon_0,1-\epsilon_0]$ for all 
 $x\in \barO$ and some $\epsilon_0\in(0,1)$, there are no 
 complications due to the boundary layer $t\searrow 0$,
 and so $u$ can be seen as 
 a ``classical'' solution over the whole $(0,T)$. This can be deduced 
 simply by using the classical Gronwall lemma (instead of the 
 uniform Gronwall lemma) in the a-priori estimates. 
 As a consequence, estimates \eqref{stz1}, \eqref{st72} and 
 \eqref{st72c} below hold in fact with $\tau=0$ in this case.
\eddos 
\beos\label{mirzel}
 It is worth stressing that the separation property~\eqref{separs}
 holding for our model is instead an open issue in the frame of the 
 standard Cahn-Hilliard system with logarithmic nonlinearity
 $f(u)$, at least in the three-dimensional setting 
 (cf.~\cite{MZ} for further remarks). 
\eddos


\section{Local strong solutions}
\label{sec:strong}

%
%


In this section, we will prove existence 
of at least one {\sl local in time}\/ classical
solution to Problem~(P). With this aim, we first 
introduce a regularization of the 
initial datum $u_0$. This is the object of the following
\bele\label{lemma:u0}
 Let $u_0$ satisfy~\eqref{EE0} and~\eqref{regou02}.
 Then, for any $\delta\in(0,1/6)$, there exists 
 $u\zzd$ such that $u\zdd \in C^{0,a}(\barO)$ 
 for any $a\in(0,1/2)$. Moreover, 
 \begin{equation} \label{uzzd0}
   -1+3\delta \le u\zzd(x) \le 1-3\delta \quad \perogni
   x\in \barO.
 \end{equation}
 Finally, we have that 
 \begin{align}\label{uzzd1} 
   & u\zzd \to u_0 \quext{weakly in }\,V,\\
  \label{uzzd2} 
   & \calE(u\zzd) \le c \big( 1 + \calE(u_0) \big) \quad \perogni \delta\in(0,1/6).
 \end{align}
\enle
\begin{proof}
First of all, we set 
\begin{equation} \label{uzzd3}
  u\zzd^{(1)}:= \min\big\{ 1 - 3\delta, \max\{ u_0, -1 + 3\delta \} \big\},
\end{equation}
so that $u\zzd^{(1)}$ satisfies a.e.~in $\Omega$ the equivalent of 
\eqref{uzzd0}. Moreover, it is clear that 
$\calE(u\zzd^{(1)})\le \calE(u_0)$ for all $\delta$.
Next, we define $z\zzd^{(1)}:=\phi(u\zzd^{(1)})=\arcsin u\zzd^{(1)}$.
Then, we proceed by singular perturbation, defining $z\zzd$ as the 
unique solution of the elliptic problem
\begin{equation} \label{uzzd4}
  z\zzd + \delta A z\zzd = z\zzd^{(1)}.
\end{equation}
Being $z\zzd^{(1)}\in V$, then, by elliptic regularity
(recall that $\Omega$ is a smooth domain),
$z\zzd\in H^3(\Omega)$ for all $\delta$. Setting 
$u\zzd:= \sin z\zzd$, a direct check permits to verify
that (at least) $u\zzd\in H^2(\Omega)$ for all $\delta$.
Hence, $u\zzd$ is H\"older continuous, as desired.
Moreover, by monotonicity of $\phi$ and a standard 
maximum principle argument, it is clear that 
\eqref{uzzd0} holds. 

The key step consists in proving \eqref{uzzd2}. 
Actually, it is obvious that $\| F(u\zzd) \|_{L^1(\Omega)} \le c$
(cf.~\eqref{regou0}).
To control the gradient term of $\calE$,
we test \eqref{uzzd4} by $A z\zzd$. We obtain
\begin{align} \no
  \frac12 \io a(u\zzd) |\nabla u\zzd|^2 &
   = \| \nabla z\zzd \|^2 
   \le \| \nabla z\zzd^{(1)} \|^2\\
 \label{uzzd5}
  & = \frac12 \io a(u\zzd^{(1)}) |\nabla u\zzd^{(1)}|^2
   \le \frac12 \io a(u_0) |\nabla u_0|^2 
   \le \calE(u_0) + c,
\end{align}
as desired. Hence, we have \eqref{uzzd2}. To conclude, we
have to prove \eqref{uzzd1}, which is however an immediate
consequence of standard weak compactness arguments.
%
%
The proof is complete.
\end{proof}
\noindent%
As a next step, we also provide a modification of 
the function $f$ given by~\eqref{f}. 
Namely, for all $\delta\in(0,1/6)$,
we take $f_\delta\in C^1(-1+\delta,1-\delta)\to \RR$ such 
that $f_\delta$ is monotone, $|f_\delta(r)| \ge |f(r)|$ 
for all $r\in(-1+\delta,1-\delta)$, and
\begin{align} \label{defifd}
  & f_\delta(r) = f(r) \quext{if }\,r\in[-1+2\delta,1-2\delta],\\
 \label{defifd2}
  & \lim_{|r|\to 1-\delta} f_\delta(r)\sign r = +\infty. 
\end{align}
It is obvious that, for any $\delta\in(0,1/6)$, a function
$f\dd$ with the above properties exists.

Finally, we modify $a$ by taking $a\dd\in C^2(\RR;\RR)$ such 
that 
\begin{equation} \label{defiad0}
  a\dd(r) = a(r) \quad\perogni r\in[-1+\delta, 1-\delta].
\end{equation}
In particular,
\begin{equation} \label{defiad}
  a\dd''(r)\ge 0, \quad \Big(\frac1{a\dd}\Big)''(r) = -1
   \quad \perogni r\in [-1+\delta,1-\delta].
\end{equation}
Outside $(-1,1)$, $a\dd$ is taken as a constant
$K_\delta>0$ (exploding as $\delta\searrow 0$),
whereas for $|r|\in (1-\delta,1)$, $a\dd$ is chosen 
in such a way to have
\begin{equation} \label{defiad2}
  1 \le a\dd(r) \le K_\delta 
   \quad \perogni r \in \RR.
\end{equation}
Then, for $\delta\in(0,1/6)$, we can consider the system
\begin{align}\label{CH1d}
  & u_{\delta,t} - \Delta w\dd = 0,
   \quext{in }\,(0,T)\times\Omega,\\
 \label{CH2d}
  & w\dd = - a\dd(u\dd) \Delta u\dd - \frac{a\dd'(u\dd)}2 |\nabla u\dd|^2 
   + f\dd(u\dd) - \lambda u\dd + \epsi u_{\delta,t},
   \quext{in }\,(0,T)\times\Omega,\\
 \label{inizd}
  & u\dd|_{t=0} = u\zzd,
  \quext{in }\,\Omega,\\
 \label{neumd}
  & \dn u\dd = \dn w\dd = 0,
   \quext{on }\, \Gamma.
\end{align}
We then have:
\bete\label{teo:loc:delta}
 Let\/ $f$ and\/ $a$ be given respectively by\/ \eqref{f}, 
 \eqref{defia},
 and let\/ $u_0$ satisfy~\eqref{EE0} and~\eqref{regou02}.
 For $\delta\in(0,1/6)$, let $u\zzd$ be defined by\/
 {\rm Lemma~\ref{lemma:u0}} and $f\dd$, $a\dd$ be given by\/
 \eqref{defifd}-\eqref{defiad2}. 
 Then, there exists one and only one 
 solution $(u\dd,w\dd)$ to system \eqref{CH1d}-\eqref{neumd}
 with the regularity
 \begin{align}\label{regoudd}
   & u\dd \in W^{1,\infty}(0,T;V') \cap H^1(0,T;V) \cap L^\infty(0,T;H^2(\Omega)), 
    \quad \epsi u\dd \in W^{1,\infty}(0,T;H),\\
  \label{regowdd}
   & w\dd \in L^\infty(0,T;V)\cap L^2(0,T;H^3(\Omega)), \quad
    \epsi w\dd \in L^\infty(0,T;H^2(\Omega)).
 \end{align}
\ente
\begin{proof}
We claim that this result is essentially a consequence of the
results of \cite{SP11}. Actually, we see that, for any 
$\delta\in (0,1/6)$, $a\dd$ satisfies the assumptions
\cite[(2.1)-(2.2)]{SP11} (where $-1$ and $1$ replace
$-2$ and $2$ in \cite[(2.2)]{SP11}) and 
\cite[(6.1)]{SP11} (where $[-1,1]$ is replaced by
$[-1+\delta,1-\delta]$, cf.~\eqref{defiad}).
 
Moreover, $f\dd$ satisfies \cite[(2.3)-(2.4)]{SP11},
with $(-1+\delta,1-\delta)$ replacing $(-1,1)$. 
Thus, we can apply \cite[Theorems~5.1, 6.1, 6.2]{SP11}
which give the existence and uniqueness of a weak solution 
$(u\dd,w\dd)$ to \eqref{CH1d}-\eqref{neumd}. The regularity
of this solution is specified 
by~\cite[(4.3)-(4.4) and (6.14)]{SP11}. More precisely,
since the initial datum $u\zzd$ is smooth and 
separated in the uniform norm from the singular values
$\pm(1-\delta)$ of $f\dd$ due to \eqref{uzzd0}, 
we have here that \cite[(6.14)]{SP11} 
holds starting from the initial
time, i.e., with $\tau=0$. Moreover,
a closer inspection of \cite[Proof of Theorem~6.2]{SP11}
(see in particular estimates~(6.19)-(6.20) therein) 
permits to see that the additional regularity 
for $u_{\delta,t}$ stated in \eqref{regoudd} holds
over $(0,T)$.
%
%
The $L^2(0,T;H^3(\Omega))$ regularity of $w\dd$ follows
from the $L^2(0,T;V)$ regularity of $u_{\delta,t}$ and elliptic
regularity estimates applied to \eqref{CH1d}.

Hence, collecting all the information coming from
the results of \cite{SP11}, we obtain exactly
\eqref{regoudd}-\eqref{regowdd}. This concludes the
proof.
\end{proof}
\noindent%
Notice now that, as a consequence of \eqref{regoudd}-\eqref{regowdd}
and of the arguments in \cite{SP11}, we have, more precisely,
the a-priori estimate
\begin{equation} \label{std11}
  \| u\dd \|_{H^1(0,T;V')} + \| u\dd \|_{L^\infty(0,T;H^2(\Omega))} 
   \le Q ( \| u\zzd \|_{H^2(\Omega)}, \delta^{-1}).
\end{equation}
%
By interpolation and embedding properties of Sobolev spaces, 
we then obtain
\begin{equation} \label{std12}
  \| u\dd \|_{C^{0,b}([0,T]\times \barO)}
   \le Q ( \| u\zzd \|_{H^2(\Omega)}, \delta^{-1}),
   \quext{for some }\,b>0.
\end{equation}
Thus, $u\dd$ is globally H\"older continuous. 
Since $u\zzd$ satisfies
 \eqref{uzzd0}, we can deduce that,
once the initial datum $u_0$ is assigned, then
for all $\delta\in(0,1/6)$ there exists 
a (computable) time $T_\delta>0$ 
depending on $\delta$ and $u_0$ such that
\begin{equation} \label{std13}
  -1 + 2\delta \le u\dd(t,x) \le 1 - 2\delta 
   \quext{for all }\,(t,x) \in [0,T_\delta]\times \barO.
\end{equation}
However, over $[-1 + 2\delta,1 - 2\delta]$, $a\dd$ coincides
with $a$ by \eqref{defiad0} and $f\dd$ coincides with $f$
by \eqref{defifd}. Hence, we have obtained the
\beco\label{cor:delta}
 Let\/ $f$ and\/ $a$ be given by\/ \eqref{f}, \eqref{defia},
 and let\/ $u_0$ satisfy~\eqref{EE0} and~\eqref{regou02}.
 Let $\delta\in(0,1/6)$ and let $u\zzd$ be given by\/
 {\rm Lemma~\ref{lemma:u0}}.
 Then, there exist a time $T_\delta>0$ depending on $u_0$ and $\delta$,  
 and a ``classical'' solution $(u\dd,w\dd)$ to~{\rm Problem~(P)},
 with initial datum $u\zzd$, over the time interval $(0,T\dd)$.
\enco


\section{A priori estimates}
\label{sec:apriori}

In this section, we derive a number of a-priori estimates
for the solutions of our system. We stress that
the procedure leading to these estimates can be rigorously
justified at least for ``classical'' solutions 
to Problem~(P) in the sense of Definition~\ref{def:classsol}.
Indeed, owing to~Remark~\ref{morereg}, these solutions can
be thought to be as smooth as we need (possibly paying the 
price of additionally regularizing $u_0$).

In particular, the estimates proved below will
hold for local strong solutions given by
Theorem~\ref{teo:loc:delta}. On the other hand, for weaker
notions of solutions the procedure below may just 
have a formal character due to insufficient regularity
of test functions. We will clarify this point
in Section~\ref{sec:weak} below.
 Here, we will proceed
assuming that everything is regular enough for our purposes.
Moreover, since the estimates we are going to derive
will have a global-in time character, with some abuse
of notation we will directly work on the time interval
$[0,T]$. The underlying extension argument will be 
also detailed in the next section.

We can now start detailing our estimates.

\smallskip

\noindent%
{\bf Energy estimate.}~~%
We test \eqref{CH1} by $w$ and \eqref{CH2} by $u_t$.
This gives rise to the energy equality \eqref{energyineq},
whence we obtain the estimate
\begin{align}\label{st11}
  & \| u_t \|_{L^2(0,T;V')}^2
   + \epsi \| u_t \|_{L^2(0,T;H)}^2
   + \| \nabla w \|_{L^2(0,T;H)}^2 \le Q(\EE_0),\\
 \label{st12}
  & \| F(u) \|_{L^\infty(0,T;L^1(\Omega))} 
   \le Q(\EE_0).
\end{align}
Thus, by definition \eqref{logpot} of $F$,
we infer in particular that
\begin{equation}\label{st13}
  -1 \le u \le 1 \quext{a.e.~in }\, 
   (0,T)\times\Omega.
\end{equation}
Moreover, recalling also \eqref{defiz}, we obtain
\begin{equation}\label{st14}
  \| u \|_{L^\infty(0,T;V)} 
   + \| z \|_{L^\infty(0,T;V)}  \le Q(\EE_0).
\end{equation}

\smallskip

\noindent%
{\bf Estimate on time derivatives.}~~%
Following the lines of \cite[Theorem~6.1]{SP11}, we 
indicate by $J$ the gradient part of the energy, i.e.,
\begin{equation} \label{defiJ}
  J:V\to [0,+\infty), \qquad
   J(u):=\io \frac{a(u)}2 | \nabla u |^2.
\end{equation}  
Then, we can (formally) compute the first 
derivative of $J$, given by 
\begin{equation} \label{Jprime}
  \duavg{J'(u),v} =
    \io \Big( a(u) \nabla u \cdot \nabla v
    + \frac{a'(u)}2 |\nabla u|^2 v \Big),   
\end{equation}  
as well as the second derivative
\begin{equation}\label{Jsecond}
  \duavg{J''(u)v,z} = 
  \io \Big( \frac{a''(u)|\nabla u|^2 vz}2
   + a'(u) v \nabla u\cdot\nabla z 
   + a'(u) z \nabla u\cdot\nabla v 
   + a(u) \nabla v\cdot\nabla z \Big).
\end{equation}  
To be more precise, if $u$ is a smooth solution 
(and in particular it is separated from singularities,
i.e., it satisfies \eqref{separs}),
formulas \eqref{Jprime} and \eqref{Jsecond}
make sense for $u = u(t)$ at any time $t\in[0,T]$
and, indeed, we have that $J'(u)\in V'$ 
and $J''(u)\in \calL(V,V')$.

From \eqref{Jsecond}, we then obtain in particular
\begin{align} \no
  \duavg{J''(u)v,v} 
  & = \io \Big( \frac{a''(u)|\nabla u|^2 v^2}2
    + 2 a'(u) v \nabla u\cdot\nabla v 
    + a(u) | \nabla v |^2 \Big)\\
 \label{Jvv}    
  & \ge \io \Big( a(u) - \frac{2 a'(u)^2}{a''(u)} \Big) 
   | \nabla v|^2.
\end{align}  
We can now test \eqref{CH1} by $t w_t$ and add the time 
derivative of \eqref{CH2} tested by $t u_t$. 
This leads to
\begin{equation}\label{Irena11}
  \frac{t}{2} \ddt \| \nabla w \|^2 
   + \frac{t \epsi}{2} \ddt \| u_t \|^2
   + t \duavg{J''(u)u_t,u_t} 
   + t \io f^{\prime}(u) u_t^2 = \lambda t \| u_t \|^2.
\end{equation}  
Hence, using \eqref{Jvv}, noting that 
\begin{equation}\label{co21}
  a(u)-\frac{2a'(u)^2}{a''(u)}
   = \frac{2}{1+3u^2},
\end{equation}
and using \eqref{st13}, 
we arrive at the inequality
\begin{align} \no
  & \ddt\Big(
   \frac{t}2 \| \nabla w \|^2
    + \frac{\epsi t}2 \| u_t \|^2 \Big)
   + t \io \Big( \frac{2}{1+3u^2} 
   | \nabla u_t |^2 \Big)
   + t \io f'(u) u_t^2\\
  \label{new11}
  & \mbox{}~~~~~
   \le \lambda t \| u_t \|^2 
    + \frac12 \| \nabla w \|^2
    + \frac{\epsi}2 \| u_t \|^2
  \le \frac{t}4 \| \nabla u_t \|^2
    + c (1+t) \| \nabla w \|^2
    + \frac{\epsi}2 \| u_t \|^2.
\end{align}  
Recalling \eqref{st11} and 
noting that $2/(1+3u^2)\ge 1/2$ since
$u$ takes values in $[-1,1]$, 
integration in time of \eqref{new11} gives,
for any $\tau \in (0,T)$,
\begin{align}\label{st21}
  & \| u_t \|_{L^\infty(\tau,T;V')} 
   + \| \nabla w \|_{L^\infty(\tau,T;H)} \le Q(\EE_0,\tau^{-1}),\\
 \label{st22}
  & \epsi \| u_t \|_{L^\infty(\tau,T;H)}
    + \| u_t \|_{L^2(\tau,T;V)} \le Q(\EE_0,\tau^{-1}).
\end{align}
%
%
%
%
%

\smallskip

\noindent%
{\bf Estimate of $f(u)$.}~~%
We test \eqref{CH2} by $u-m$. Integrating
by parts the terms depending on $a$, we obtain
\begin{equation}\label{conto51}
  \io \Big( a(u) + \frac{a'(u)}2 ( u - m ) \Big) | \nabla u |^2
   + \io f(u) ( u - m )
    =\big( w + \lambda u - \epsi u_t, u - m \big).
\end{equation}
We have to estimate some terms. Firstly,
proceeding as in~\cite[Appendix]{MZ}, it 
is not difficult to prove that
\begin{equation}\label{conto52}
  \io f(u) ( u - m )
   \ge \frac12 \| f(u) \|_{L^1(\Omega)} - c_m.
\end{equation}
We notice that assumption~\eqref{regou02} is used here.

We now observe that
\begin{equation}\label{co31}
  \io \Big( a(u) + \frac{a'(u)}2 ( u - m ) \Big) | \nabla u |^2
  = \io \Big( \frac{2u(u-m)}{(1-u^2)^2} + \frac2{1-u^2} \Big)
   |\nabla u|^2.
\end{equation}
Moreover, there exist constants $\kappa_m>0$, $c_m\ge 0$ 
such that
\begin{equation}\label{co32}
  \frac{2u(u-m)}{(1-u^2)^2} 
   \ge \kappa_m \frac{1}{(1-u^2)^2}
   - c_m
   \ge \kappa_m \frac{u^2}{(1-u^2)^2}
   - c_m.
\end{equation}
Thus, recalling \eqref{defia}, we get 
\begin{equation}\label{co31x}
  \io \Big( a(u) + \frac{a'(u)}2 ( u - m ) \Big) | \nabla u |^2
   \ge \kappa_m \| \nabla f(u) \|^2 - c_m \| \nabla u \|^2.
\end{equation}
Finally, noting that $u$ takes values in $[-1,1]$,
using estimate \eqref{st14}, and applying the 
Poincar\'e-Wirtinger inequality, we have
\begin{equation}\label{conto54}
  \big( w + \lambda u - \epsi u_t, u - m \big)
   = \big( w - w\OO + \lambda u - \epsi u_t, u - m \big) 
    \le c \big( \| \nabla w \| + 1 + \epsi \| u_t \| \big).
\end{equation}
Thus, collecting the above considerations, \eqref{conto51} gives
\begin{equation}\label{conto51b}
  \| \nabla f(u) \|^2
   + \| f(u) \|_{L^1(\Omega)}
    \le c \big( \| \nabla u \|^2 + \| \nabla w \| 
     + 1 + \epsi \| u_t \| \big).
\end{equation}
Squaring \eqref{conto51b}, using
\eqref{conto52}-\eqref{conto54}, and integrating
in time, we arrive at
\begin{equation}\label{st31}
  \| \nabla f(u) \|_{L^4(0,T;H)} 
   + \| f(u) \|_{L^2(0,T;L^1(\Omega))}
   \le Q(\EE_0),
\end{equation}
whence, more precisely,
\begin{equation}\label{st32}
  \| f(u) \|_{L^2(0,T;V)} \le Q(\EE_0).
\end{equation}
Taking instead the essential supremum of~\eqref{conto51b}
as $t$ ranges in $(\tau,T)$, and using \eqref{st21}-\eqref{st22},
a straighforward modification of the above procedure 
leads to 
\begin{equation}\label{st35}
  \| f(u) \|_{L^\infty(\tau,T;V)} \le Q(\EE_0,\tau^{-1})
   \quad \perogni \tau\in(0,T).
\end{equation}
Next, integrating \eqref{CH2} in space, 
using \eqref{aprimo}, \eqref{st13}, \eqref{defia},
and noting that $(u_t)\OO\equiv 0$, we get
\begin{align}\no
  | w\OO |
   & \le c \Big( \io \frac{a'(u)}2 | \nabla u |^2 
    + \| f(u) \|_{L^1(\Omega)}
     + 1  \Big)\\
 \label{conto51c}    
  & \le c \Big( \| \nabla f(u) \|^2 
    + \| f(u) \|_{L^1(\Omega)}
     + 1   \Big).
\end{align}
Thus, squaring, integrating in time, using \eqref{st31}, 
and recalling the last~\eqref{st11}, we infer
\begin{equation}\label{st34}
  \| w \|_{L^2(0,T;V)} \le Q(\EE_0).
\end{equation}
Taking the essential supremum in \eqref{conto51c} as $t\in(\tau,T)$,
and recalling \eqref{st21} and \eqref{st35}, we also get
\begin{equation}\label{st36}
  \| w \|_{L^\infty(\tau,T;V)} \le Q(\EE_0,\tau^{-1})
   \quad \perogni \tau\in(0,T).
\end{equation}

\smallskip

\noindent%
{\bf Estimate of $z$.~~}%
We consider the equivalent formulation
\eqref{CH2z} and test it by $-\Delta z$. This gives
\begin{equation}\label{coz1}
  \io 2 \phi'(u) | \Delta z |^2
   + \big( f'(u) \nabla u,\nabla z \big)
   = \big( w + \lambda u - \epsi u_t, - \Delta z \big).
\end{equation}
Thus, using the monotonicity of $f$ and $\phi$, 
noting that $\phi'(u)\ge 1/2$ for all $u\in(-1,1)$,
and that
$\big( f'(u) \nabla u,\nabla z \big) \geq 2 \| \nabla z \|^2 $,
we can control the \rhs\ this way:
\begin{equation}\label{coz2}
  \big( w + \lambda u - \epsi u_t, - \Delta z \big)
   \le \frac12 \| \Delta z \|^2
    + c \big( \| w \|^2 + 1 + \epsi^2 \| u_t \|^2 \big).
\end{equation}
Hence, integrating \eqref{coz1} in time
and recalling~\eqref{st11}, \eqref{st14}
and \eqref{st34}, we arrive at
\begin{equation}\label{stz2}
  \| z \|_{L^2(0,T;H^2(\Omega))} \le Q(\EE_0).
\end{equation}
Taking instead the (essential) supremum of \eqref{coz1} as $t$ 
ranges in $(\tau,T)$ for $\tau>0$ and using 
\eqref{st21}-\eqref{st22} and \eqref{st36}, we obtain
\begin{equation}\label{stz1}
  \| z \|_{L^\infty(\tau,T;H^2(\Omega))} \le Q(\EE_0,\tau^{-1}).
\end{equation}
The above relations permit to improve also the bounds on
$u$. Actually, computing directly the Laplacean
of $u = \sin z$ and using \eqref{st13} together with
the Gagliardo-Nirenberg 
inequality (cf., e.g., \cite[Theorem p.~125]{Ni})
\begin{equation}\label{ineq:gn}
  \| \nabla y \|_{L^4(\Omega)} \le c\OO \| y \|_{H^2(\Omega)}^{1/2}
   \| y \|_{L^\infty(\Omega)}^{1/2}
   + \| y \| \qquad \perogni y \in H^2(\Omega),
\end{equation}
it is not difficult to arrive at 
\begin{align}\label{stnewu1}
  & \| u \|_{L^2(0,T;H^2(\Omega))} \le Q(\EE_0),\\
 \label{stnewu2}
  & \| u \|_{L^\infty(\tau,T;H^2(\Omega))} \le Q(\EE_0,\tau^{-1}).
\end{align}

\smallskip

\noindent%
{\bf First entropy estimate and separation property.}~~%
The estimates obtained up to this moment yield a control 
of the functions $z$ and $u$ up to their second space derivatives
(cf.~\eqref{stz2}-\eqref{stz1} 
and \eqref{stnewu1}-\eqref{stnewu2}), and of the nonlinear
term $f(u)$ up to its first space derivatives
(cf.~\eqref{st32} and \eqref{st35}). However,
this still seems not sufficient to pass to the limit
in the equation \eqref{CH2}, even if its equivalent
formulation \eqref{CH2z} is considered. Indeed, 
from \eqref{st32} and \eqref{st35} we get 
a control of the term $f(u)$, that explodes logarithmically
fast as $|u|\nearrow 1$. On the other hand,
even in  formulation \eqref{CH2z} one faces the term 
$\phi'(u)$ which is much more singular since it
explodes as a negative power of $(1-u^2)$. 
To control it in some $L^p$-norm we need more refined 
estimates of the so-called {\sl entropy}\/ type and, 
in particular, we need to refer to the formulation 
\eqref{CH2v} (we recall that all formulations 
are equivalent, at least for sufficiently smooth solutions).
Usage of this technique requires the convexity assumption 
on $\Omega$ asked in the statement of Theorem~\ref{teo:main}.

The basic tool we need consists in an
integration by parts formula due to Dal Passo, 
Garcke and Gr\"un (\cite[Lemma 2.3]{DpGG}):
\bele\label{lemma:dpgg}
 Let $h\in W^{2,\infty}(\RR)$ and $z\in W$. Then,
 \begin{align} \no
   & \io h'(z) |\nabla z|^2 \Delta z
    = -\frac13 \io h''(z) |\nabla z|^4\\
  \label{byparts}
   & \mbox{}~~~~~
    + \frac23 \io h(z) \big( |D^2 z|^2 - | \Delta z|^2 \big)
    + \frac23 \iga h(z) II( \nabla z ),
 \end{align}  
 where $II(\cdot)$ denotes the second fundamental form
 of $\Gamma$. 
\enle
\noindent%
Then, we test \eqref{CH2v} by 
\begin{equation}\label{co63}
  - \Delta\mciapo(v) = -\dive (m(v) \nabla v)
   = - m(v) \Delta v - m'(v) | \nabla v|^2,
\end{equation}
with the function $m$ to be chosen later. This gives
\begin{align} \no
  & \io m(v) | \nabla v |^2 
   + \io \bigg(  m(v) |\Delta v|^2 
   - \frac{m'(v)j(v)}2 |\nabla v|^4
   + \Big( m'(v) - \frac{m(v)j(v)}2 \Big) 
            \Delta v | \nabla v|^2 \bigg)\\
 \label{co64}
  & \mbox{}~~~~~
   = \big( m(v) \nabla v, \nabla w + \lambda \nabla u \big)
    + \big( \epsi u_t, m(v) \Delta v + m'(v) | \nabla v|^2 \big).
\end{align}  
Now, applying Lemma~\ref{lemma:dpgg} to
the last integral on the \lhs\ of \eqref{co64} we infer
\begin{align}\no
  & \io \Big( m'(v) - \frac{m(v)j(v)}2 \Big) \Delta v | \nabla v|^2
   = - \frac13 \io \Big( m''(v) - \frac{m'(v)j(v)}2 - \frac{m(v)j'(v)}2 \Big) | \nabla v|^4\\
 \label{co65}
  & + \frac23 \io \Big( m(v) - \frac{\mjciapo(v)}2 + K \Big)\big( | D^2 v|^2 - | \Delta v |^2\big)
   + \frac23 \iga \Big( m(v) - \frac{\mjciapo(v)}2 + K \Big) II(\nabla v),
\end{align}
where $K>0$ is an integration constant that will be chosen
later on, and the notation \eqref{deficiapo} is used.

Substituting \eqref{co65} into \eqref{co64}, we get on 
the \lhs\ the following ``hopefully good''
terms
\begin{equation}\label{co66}
  \io m(v) |\nabla v|^2 + \io m(v) |\Delta v|^2 
   + \frac13 \io \Big( - m''(v) - m'(v) j(v) + \frac{m(v)j'(v)}2 \Big) | \nabla v|^4,
\end{equation}
where we notice that 
\begin{equation}\label{co67}
  j'(v) = \frac{2e^v}{(e^v+1)^2}.
\end{equation}
We can now specify our choice of $m$ as
\begin{equation}\label{defim}
  m(v)=\frac1{2(1+v^2)}.
\end{equation}
Actually, this expression arises since we need $m$ to
decay not too fast at infinity (otherwise we 
do not get enough information from it), but at the same time 
we need it to be summable (cf.~\eqref{co75} below).
The above choice gives
\begin{equation}\label{mprim}
  \mciapo(v) =\frac12 \arctan v,
\end{equation}
as well as
\begin{equation}\label{mder}
  m'(v)=-\frac{v}{(1+v^2)^2}, \qquad
   m''(v)=\frac{3v^2-1}{(1+v^2)^3}.
\end{equation}
Now, noting that $mj'\ge 0$, we can easily observe that
\begin{equation}\label{co68}
  - m''(v) - m'(v) j(v) + \frac{m(v)j'(v)}2
  \ge \frac{(v + v^3)j(v) + 1 - 3 v^2}{(1+v^2)^3}.
\end{equation}
Now, recalling \eqref{co61}, for a suitable $M>0$ 
(e.g., we can take $M=12$ here), a direct computation shows 
that
\begin{equation}\label{co68bis}
  \frac{(v + v^3)j(v) + 1 - 3 v^2}{(1+v^2)^3}
   \ge \frac14 \frac{1 + |v|^3}{(1+v^2)^3}
   \quad\perogni |v|\ge M.  
\end{equation}
On the other hand,
\begin{equation}\label{co68ter}
  \frac{(v + v^3)j(v) + 1 - 3 v^2}{(1+v^2)^3}
   \ge \frac1{(1+v^2)^3} - \frac{3 v^2}{(1+v^2)^3}
   \ge \kappa_M \frac{1 + |v|^3}{(1+v^2)^3} 
    - c_M \frac{|v|}{(1+v^2)^3}
   \quad\perogni |v| \le M.  
\end{equation}
Summarizing, we have 
\begin{equation}\label{co68fine}
  - m''(v) - m'(v) j(v) + \frac{m(v)j'(v)}2
   \ge \kappa \frac{1+|v|^3}{(1+v^2)^3}
    - c \frac{|v|}{(1+v^2)^3} \chi_{\{|v|\le M\}},
\end{equation}
where $\chi$ denotes the characteristic function.

Now, we can choose $K$ so large that the function 
\begin{equation}\label{co75}
  \Big( m(v) - \frac{\mjciapo(v)}2 + K \Big)
\end{equation}
is strictly positive (and bounded, of course). 
Then, thanks the {\sl convexity}\/ assumption on $\Omega$
(that entails positive definiteness of the second
fundamental form), the latter two terms in \eqref{co65} are 
positive. 

Collecting these observations, we can deduce from
\eqref{co64} the estimate
\begin{equation}\label{co64b}
  \io m(v) | \nabla v |^2  
   + \io m(v) |\Delta v|^2 
  + \kappa \io \frac{1+|v|^3}{(1+v^2)^3} | \nabla v|^4
  \le \sum_{i=1}^4 T_i,
\end{equation}
and we have to control the ``bad'' terms on the \rhs:
\begin{align}\label{co72}
  & T_1:= c \int_{\{|v|\le M\}} \frac{|v|}{(1+v^2)^3} | \nabla v|^4,\\
 \label{co73}
  & T_2:= \io m(v) \nabla v \cdot \nabla w,\\
 \label{co74}
  & T_3:= \lambda \io m(v) \nabla v \cdot \nabla u,\\
 \label{co74x}
  & T_4:= \big( \epsi u_t, m(v) \Delta v + m'(v) | \nabla v|^2 \big).
\end{align}
To do this, we first notice that
\begin{equation}\label{co74b}
  | T_2 + T_3 | 
   \le \frac14 \io m(v) | \nabla v |^2
    + c \big( \| \nabla w \|^2 + \| \nabla u \|^2 \big).
\end{equation}
Next, we have to control $T_1$. Then, we can note
that there exists $\delta\in(0,1)$, depending only on $M$,
such that the restriction
$f:[-1+\delta,1-\delta]\to[-M,M]$ is bijective
and Lipschitz continuous together with its 
inverse $j$. 
%
%
Thus, using \eqref{st13} and
inequality \eqref{ineq:gn}, we deduce
%
%
%
%
\begin{align}\no
  | T_1 | & \le c \int_{\{|v|\le M\}} | \nabla v |^4
    \le c \int_{\{|v|\le M\}} f'(u)^4 | \nabla u |^4\\
 \no
   & \le c \int_{\{|u|\le 1-\delta\}} \frac1{(1-u^2)^4} | \nabla u|^4
     \le c_\delta \int_{\{|u|\le 1-\delta\}} | \nabla u |^4\\
 \label{co72b}
  & \le c_\delta \io | \nabla u |^4
   \le c_\delta \| u \|_{L^\infty(\Omega)}^2
     \big( 1 + \| u \|_{H^2(\Omega)}^2 \big)
     \le c_\delta \big( 1 + \| u \|_{H^2(\Omega)}^2 \big).
\end{align}
Finally, we have to control $T_4$. We have
\begin{equation}\label{co74d}
  | T_4 |  \le \frac14 \io m(v) | \Delta v |^2
    + c_\sigma \epsi^2 \| u_t \|^2
    + \sigma \io m'(v)^2 |\nabla v|^4,
\end{equation}
for small $\sigma>0$ to be chosen later. Then,
we can go on as follows:
\begin{equation}\label{co74d2}
   \sigma \io m'(v)^2 |\nabla v|^4
    = \sigma \io \frac{ v^2 }{ (1 + v^2)^4 }|\nabla v|^4
    \le \sigma \io \frac{ 1 + |v|^3 }{ (1 + v^2)^3 }|\nabla v|^4.
\end{equation}
Hence, for $\sigma$ sufficiently small, the last
integral is controlled by the last term on
the \lhs\ of \eqref{co64b}. On the other hand, the
first term on the \rhs\ of \eqref{co74d} is
controlled by the corresponding one on the
\lhs\ of \eqref{co64b}.

Now, a direct computation permits to see that
\begin{equation}\label{co71}
  \io \frac{1+|v|^3}{(1+v^2)^3} | \nabla v|^4
  \ge \io \frac{|v|^3}{(1 + v^2)^3} \frac{|v|}{(1+ v^2)^{1/2}} | \nabla v|^4
   \ge \kappa \io \big| \nabla (1+v^2)^{1/8} \big|^4.
\end{equation}
Thus, collecting \eqref{co74b}-\eqref{co71},
\eqref{co64b} gives
\begin{align}\no
  & \io m(v) | \nabla v |^2 
   + \io m(v) | \Delta v |^2 
   + \io \big| \nabla (1+v^2)^{1/8} \big|^4\\
 \label{co64b2}  
  & \mbox{}~~~~~
   \le c \big( \| \nabla w \|^2 + \| \nabla u \|^2 
    + 1 + \| u \|_{H^2(\Omega)}^2 
    + \epsi^2 \| u_t \|^2 \big).
\end{align}
Hence, integrating in time, and using \eqref{st11}, \eqref{st14}
and \eqref{stnewu1}, we obtain
\begin{equation}\label{st71}
  \big\| m^{1/2}(v) \Delta v \big\|_{\LDH} 
   + \big| (1+v^2)^{1/8} \big|_{L^4(0,T;W^{1,4}(\Omega))} 
  \le Q(\EE_0).
\end{equation}
On the other hand, taking the essential supremum
in \eqref{co64b2} as $t$ ranges in $(\tau,T)$ for
$\tau>0$, and using \eqref{st21}, \eqref{st22} and
\eqref{stnewu2}, we arrive at
\begin{equation}\label{st72}
  \big\| m^{1/2}(v) \Delta v \big\|_{L^\infty(\tau,T;H)}
   + \big| (1+v^2)^{1/8} \big|_{L^\infty(\tau,T;W^{1,4}(\Omega))} 
  \le Q(\EE_0,\tau^{-1}).
\end{equation}
By the continuous embedding 
$W^{1,4}(\Omega)\subset L^\infty(\Omega)$ we 
have in particular
\begin{equation}\label{st72b}
  \big| (1+v^2)^{1/8} \big|_{L^\infty((\tau,T)\times \Omega)} \le Q(\EE_0,\tau^{-1}).
\end{equation}
In terms of $u$ the above estimate gives rise to
the separation property
\begin{equation}\label{st72c}
  -1 + \epsilon \le u(x,t) \le 1 - \epsilon 
   \quext{a.e.~in }\,(\tau,T)\times \Omega,
\end{equation}
for all $\tau>0$, with $\epsilon>0$ depending on $\tau$.
%
%
%

\smallskip

\noindent%
{\bf Refined entropy estimate.}~~%
We repeat the entropy estimate of before taking now,
in place of \eqref{defim},
\begin{equation}\label{defimp}
  m(v)=\frac1{2(1+v^2)^p},
\end{equation}
where the choice of $p\in (1/2,1]$ will be 
made precise later on. Then, we have
%
%
%
\begin{equation}\label{mderp}
  m'(v)=-\frac{pv}{(1+v^2)^{p+1}}, \qquad
   m''(v)=\frac{(2p^2+p)v^2-p}{(1+v^2)^{p+2}}.
\end{equation}
Thus, a straighforward modification of 
\eqref{co68}-\eqref{co68ter} leads to
\begin{equation}\label{co68finep}
  - m''(v) - m'(v) j(v) + \frac{m(v)j'(v)}2
   \ge \kappa \frac{1+|v|^3}{(1+v^2)^{p+2}}
    - c \frac{|v|}{(1+v^2)^{p+2}} \chi_{\{|v|\le M\}},
\end{equation}
whereas the equivalent of \eqref{co71} gives now rise to
\begin{equation}\label{co71p}
  \io \frac{1+|v|^3}{(1+v^2)^{p+2}} | \nabla v|^4
   \ge \kappa \io \big| \nabla (1+v^2)^{\frac38-\frac{p}4} \big|^4.
\end{equation}
Now, let us observe that, since $p>1/2$, then the function 
$\mjciapo$ is bounded. Thus, we can still take $K$ so large,
depending of course on $p$, that the function in 
\eqref{co75} is strictly positive. 
Moreover, the terms corresponding to $T_j$, $j=1,\dots,4$, can
be controlled similarly as before. Thus, integrating
in time the $p$-analogue of \eqref{co64b}, we obtain
\begin{equation}\label{co81b}
  \iTo \big| \nabla (1+v^2)^{\frac38-\frac{p}4} \big|^4
   \le Q(\EE_0).
\end{equation}
Let us now test \eqref{CH2v} by $v$. This gives 
\begin{equation}\label{conew1}
  \| \nabla v \|^2 
   + \| v \|^2 
   + \io \frac{j(v) v}2 | \nabla v |^2
  \le \big( w + \lambda u - \epsi u_t, v \big).
\end{equation}   
Then, integrating in time,
noting that $1 + j(v)v\ge \kappa (1 + v^2)^{1/2}$, 
recalling \eqref{st11} and using H\"older's and Young's inequalities
to estimate the \rhs, we arrive at 
\begin{equation}\label{co82}
  \iTo (1+v^2)^{1/2} | \nabla v|^2 \le Q(\EE_0).
\end{equation}
Then, let us define
\begin{equation}\label{co83}
  \Omega^+(t):=\big\{x\in\Omega:~|v(x,t)|\ge 1\big\}, \qquad
   \Omega^-(t):=\big\{x\in\Omega:~|v(x,t)|\le 1\big\}.
\end{equation}
Let us now see that \eqref{co81b} and \eqref{co82} permit to
prove higher integrability properties of $\nabla v$.
Firstly, recalling \eqref{st14} and \eqref{stnewu1}, 
and using inequality \eqref{ineq:gn}, we obtain
\begin{equation}\label{co84-}
  \| u \|_{L^4(0,T;W^{1,4}(\Omega))} \le Q(\EE_0).
\end{equation}
Hence, noting that $j$ is locally Lipschitz continuous 
with its inverse, we get an analogous information for $v$ in
the space-time set where it is small:
\begin{equation}\label{co84}
  \iTT \int_{\Omega^-(t)} | \nabla v |^4 \le Q(\EE_0).
\end{equation}
On the other hand, in $\Omega^+(t)$ we can write, for $\eta>0$,
\begin{align}\no
  | \nabla v |^2 & = \big( v ( 1 + v^2 )^{-\eta} \nabla v \big) \cdot 
   \Big( \frac{( 1 + v^2 )^{\eta}} v \nabla v \Big)\\
 \label{co85} 
   & = c_\eta \nabla ( 1 + v^2 )^{1-\eta} \cdot 
   \Big( \frac{( 1 + v^2 )^{\eta}} v \nabla v \Big).
\end{align}
Thus, choosing $\eta=\frac58+\frac{p}4$, we have
$1-\eta=\frac38-\frac{p}4$. Then, we can take $p=\frac12+\epsilon$, 
with $\epsilon>0$ as small as we want. We then obtain that 
$\eta-\frac14=\frac12+\frac\epsilon4$.
Thus, recalling \eqref{co81b} and \eqref{co82}, we arrive at
\begin{equation}\label{co86}
  | \nabla v |^2 \le c_\eta \bigg| \underbrace{\nabla ( 1 + v^2 )^{1-\eta}}_{L^4}
    \cdot 
   \underbrace{\big( ( 1 + v^2 )^{1/4} \nabla v \big)}_{L^2}
   \underbrace{\frac{( 1 + v^2 )^{\eta - \frac14}}{v}}_{L^r} \bigg|,
\end{equation}
where we can take $r$ strictly greater than $4$ since 
for~$|v|\ge 1$ we have
\begin{equation}\label{co87}
  \frac{( 1 + v^2 )^{\eta - \frac14}}{|v|} 
  = \frac{( 1 + v^2 )^{\frac12+\frac\epsilon4}}{|v|} 
  \sim |v|^{\frac{\epsilon}2},
\end{equation}
and we know from \eqref{st32} that $v$ is controlled
in $L^2(0,T;H)$. Thus, estimating the
\rhs\ of \eqref{co86} by Young's inequality,
integrating first over $\Omega^+(t)$ and
then for $t\in(0,T)$, we obtain that
\begin{equation}\label{co88}
  \iTT \int_{\Omega^+(t)} | \nabla v |^q \le Q(\EE_0),
   \quext{for some }\,q~\text{{\bf strictly} larger than }\,2,
\end{equation}
and, of course, combining with \eqref{co84},
\begin{equation}\label{co88b}
  \iTo | \nabla v |^q \le Q(\EE_0),
   \quext{for some }\,q~\text{{\bf strictly} larger than }\,2.
\end{equation}
%


\section{Existence and uniqueness of weak solutions}
\label{sec:weak}

We detail here the proof of Theorem~\ref{teo:main}, which
is largely based on the estimates derived in the previous
section. As a first step, however, we show uniqueness, 
which works similarly to \cite{SP11}. Indeed, the key 
assumption \cite[(6.1)]{SP11} is satisfied by our function~$a$
(we have, indeed, that $(1/a)''\equiv -1$).

Then, let us take a couple of weak 
solutions $(u_1,w_1)$ and $(u_2,w_2)$
originating from the same initial datum $u_0$,
and assume that both are ``classical'' (and in particular
satisfy the separation property \eqref{separs}),
for strictly positive times.  
 
Setting $(u,w):=(u_1,w_1)-(u_2,w_2)$, we can write both 
\eqref{CH1w} and \eqref{CH2w} for the two solutions
and take the difference. Using notation \eqref{defiJ}
we get
\begin{align}\label{CH1di}
  & u_t + A w = 0,\\
 \label{CH2di}
  & w = J'(u_1) - J'(u_2) + f(u_1) - f(u_2) - \lambda u + \epsi u_t.
\end{align}
Then, we can test \eqref{CH1di} by
$A^{-1}u$, \eqref{CH2di} by $u$, and take the difference.
Actually, the operator $A$ is invertible as it is restricted
to $0$-mean valued functions (as in the case of $u$ due to 
conservation of mass). Noting that
\begin{equation} \label{contod0}
  \big( Aw, A^{-1} u \big)
   = \big( A (w-w\OO), A^{-1} u \big)
   = ( w - w\OO, u )
   = ( w, u ),
\end{equation}  
we then obtain
\begin{equation} \label{contod1}
  \frac12 \ddt \Big( \| u \|_{V'}^2 + \epsi \| u \|^2 \Big)
   + \duavg{J'(u_1) - J'(u_2), u}
   + \big( f(u_1) - f(u_2), u \big) 
   \le \lambda \| u \|^2.
\end{equation}  
Then, recalling \eqref{Jvv} and \eqref{co21}
and using monotonicity of $f$, we arrive~at 
\begin{equation} \label{contod2}
  \frac12 \ddt \Big( \| u \|_{V'}^2 + \epsi \| u \|^2 \Big)
   + \frac12 \| \nabla u \|^2
  \le \lambda \| u \|^2.
\end{equation}  
Noting that, by the Poincar\'e-Wirtinger inequality,
\begin{equation} \label{contod3}
  \lambda \| u \|^2  
   \le \frac14 \| \nabla u \|^2
   + c \| u \|_{V'}^2,
\end{equation}
we can integrate \eqref{contod2}
over $(\tau,T)$ for $\tau>0$. Using Gronwall's lemma, 
we obtain
\begin{equation} \label{fromSP11-2}
  \| u_1 - u_2 \|_{L^\infty(\tau,T;V')}^2
   \le c(T) \| u_1(\tau) - u_2(\tau) \|_{V'}^2,
\end{equation}   
where $c(T)$ is independent of $\tau$. Then, uniqueness
follows by taking the limit $\tau\searrow 0$ and owing
to continuity of weak solutions with values in $V'$ (which
is an obvious consequence of \eqref{regou4}). 

\smallskip

Let us now switch to existence. To start with, 
we approximate the initial datum $u_0$ as specified in
Lemma~\ref{lemma:u0}. Then, thanks to Corollary~\ref{cor:delta},
for any $\delta\in (0,1/6)$, there exists a ``classical''
solution $(u\dd,w\dd)$ to Problem~(P)
defined {\sl at least}\/ on the time interval
$(0,T\dd)$, where $T\dd$ depends on $u_0$ and $\delta$. 
Actually, in principle, we may have that $T\dd\searrow 0$ as we let 
$\delta\searrow 0$. On the other hand, the forthcoming argument
will exclude this eventuality and show that, in fact,
$(u\dd,w\dd)$ can be extended up to the final time $T$.
 
Indeed, let us denote as $T\dmax$ the maximum time up to 
which $(u\dd,w\dd)$ can be extended in the form of
a ``classical'' solution; namely,
\begin{equation}\label{defitdmax}
  T\dmax:=\sup\big\{ S\in(0,T]:~ u\dd~\text{admits a ``classical'' extension over $(0,S)$}\big\}.
\end{equation}
Due to uniqueness proved above, all extensions of 
$(u\dd,w\dd)$ can be ``glued'' together. Consequently,
there exists a (unique) maximal classical extension $(u\dmax,w\dmax)$ 
defined over $(0,T\dmax)$. We claim that $T\dmax = T$, and,
to prove this claim, we proceed as usual by contradiction.
Actually, due to \eqref{regoust}-\eqref{regowst} and \eqref{separs},
for any $S\in (0,T\dmax)$ we have
\begin{align}\label{regoumax}
  & \| u\dmax \|_{H^1(0,S;V) \cap W^{1,\infty}(0,S;V') \cap L^\infty(0,S;H^2(\Omega))}
   \le C(S), \qquad
   \epsi \| u\dmax \|_{W^{1,\infty}(0,S;H)}\le C(S),\\
 \label{regowmax}
  & \| w\dmax \|_{L^\infty(0,S;V) \cap L^2(0,S;H^3(\Omega))} \le C(S), \qquad
   \epsi \| w\dmax \|_{L^\infty(0,S;H^2(\Omega))}\le C(S),\\
 \label{separmax}
  & - 1 + \epsilon(S) \le u\dmax(t,x) \le 1 - \epsilon(S)
  \quad \perogni (t,x) \in [0,S]\times \barO,
\end{align}
where $C(S),\epsilon(S)>0$ and it may be $C(S)\nearrow\infty$ and 
$\epsilon(S)\searrow 0$ as $S\nearrow T\dmax$.
On the other hand, since $(u\dmax,w\dmax)$ is a ``classical'' 
solution, it satisfies the a-priori estimates of the previous section 
on the time interval $(0,T\dmax)$. Then, thanks to \eqref{st72c},
we have that 
\begin{equation}\label{separunif}
   - 1 + \bar\epsilon \le u\dmax(t,x) \le 1 - \bar\epsilon
  \quad \perogni (t,x) \in [\tau,S]\times \barO,
\end{equation}
and for all $0<\tau<S<T\dmax$, with $\bar\epsilon$ independent
both of $S$ and of $\delta$.
To be more precise, we have that  
\begin{equation}\label{separunif2}
  \bar\epsilon^{-1} = Q(\calE(u\zzd), \tau^{-1})
   \le Q(\EE_0, \tau^{-1}),
\end{equation}  
where the second inequality is a consequence of 
\eqref{uzzd2}.
     
Analogously, we have estimates of the norms in 
\eqref{regoumax}-\eqref{regowmax} over the time interval $(\tau,S)$
by a constant $C$ independent of $\delta$ and $S$.
Consequently, we obtain
\begin{equation}\label{extens}
  \esiste \baru = \lim_{t\nearrow T\dmax} u\dmax(t,\cdot).
\end{equation}
To be more precise, this limit is 
reached in the weak topology of $H^2(\Omega)$. Indeed, it is
a consequence of \eqref{regoumax} that 
$u\dmax\in C_w([0,T\dmax];H^2(\Omega))$.
Thus, $- 1 + \bar\epsilon \le \baru(x) \le 1 - \bar\epsilon$ for all
$x\in \barO$ and we can use $\baru$ as a new 
``initial'' datum and extend the solution $(u\dmax,w\dmax)$
beyond the time $T\dmax$. Moreover, the extension is still
a classical solution since $\baru$ is H\"older continuous 
and uniformly separated from $-1$ and $1$. This contradicts
the maximality of $T\dmax$ and of $(u\dmax,w\dmax)$.
Hence, we necessarily have that $T\dmax=T$.

To conclude the proof, we need to show that we can take 
the limit $\delta\searrow 0$ and obtain a weak solution to
Problem~(P). With this purpose we rename simply as $(u\dd,w\dd)$ 
the maximal solution obtained in the previous part 
(which is now defined in the whole $(0,T)$), and observe
that, thanks to estimates \eqref{st11}, \eqref{st14},
\eqref{stz2}, and \eqref{st32}, there hold the following 
convergence relations:
\begin{align}\label{conv11}
  & u\dd \to u \quext{weakly star in }\,\HUVp \cap \LIV \cap L^\infty((0,T)\times\Omega),\\
 \label{conv11b}
  & \epsi u\dd \to \epsi u \quext{weakly in }\,\HUH,\\
 \label{conv12}
  & w\dd \to w \quext{weakly in }\,\LDV,\\
 \label{conv13}
  & z\dd \to z \quext{weakly star in }\,\LIV \cap \LDHD,\\
 \label{conv14}
  & v\dd=f(u\dd) \to v \quext{weakly in }\,\LDV,
\end{align}
for suitable limit functions $u,w,z,v$. The above properties,
as well as the ones that will follow, are to be intended
up to the extraction of (non-relabelled) subsequences 
of $\delta\searrow 0$. We then immediately see that relation
\eqref{CH1} passes to the limit. However, since it is only
$w\in \LDV$, the Laplace operator with the boundary
condition $\dn w=0$ have to be interpreted in the weak
form through the operator $A$ (cf.~\eqref{defiA}).
  
Next, applying the Aubin-Lions lemma, \eqref{conv11} gives
\begin{equation}\label{conv21}
  u\dd \to u \quext{strongly in }\, L^p((0,T)\times\Omega),
   \quad \perogni p\in[1,\infty).
\end{equation}
Thus, a standard monotonicity argument (see, e.g., 
\cite[Prop.~1.1, p.~42]{barbu})
permits to infer from \eqref{conv14} 
that $v=f(u)$ almost everywhere. Moreover, by
the generalized Lebesgue's theorem, we get more 
precisely
\begin{equation}\label{conv22}
  v\dd=f(u\dd) \to v=f(u) \quext{strongly in }\, L^p((0,T)\times\Omega),
   \quad \perogni p\in[1,2).
\end{equation}
Now, as a consequence of estimate \eqref{co88b}, we obtain
\begin{equation}\label{conv23}
  \nabla v\dd \to \nabla v \quext{weakly in }\, L^p((0,T)\times\Omega),
   \quext{for some }\, p>2.
\end{equation}
Collecting \eqref{conv21} and \eqref{conv23}, we infer
\begin{equation}\label{conv24}
  \frac{u\dd}2 |\nabla v\dd|^2 \to \Phi \quext{weakly in }\, L^p((0,T)\times\Omega),
   \quext{for some }\, p>1,
\end{equation}
and for a suitable limit function $\Phi$. 

We can now write equation \eqref{CH2v} for the solution 
$(u\dd,w\dd)$ and see that, thanks to the above
convergence relations, all terms pass to the limit. 
Actually, what we obtain for $\delta\searrow 0$ is
\begin{equation}\label{CH2vlim}
  w = - \Delta v + v + \Phi  - \lambda u + \epsi u_t,
\end{equation}
at least in the distributional sense.
%
%
To conclude the proof, we need to identify 
the function $\Phi$. To this aim, we notice that, 
by \eqref{conv11}, \eqref{conv12}, \eqref{conv24}, and
a comparison of terms in \eqref{CH2w},
\begin{equation}\label{conv31}
  \| -\Delta v\dd \|_{L^p((0,T)\times\Omega)} \le c,
   \quext{for some }\, p>1.
\end{equation}
Thus, we get that, also in the limit, 
$-\Delta v$ (the distributional Laplacean)
lies in $L^p((0,T)\times\Omega)$ for some $p>1$.
More precisely, thanks to the no-flux condition and 
to elliptic regularity, we deduce
\begin{equation}\label{conv32}
  v\dd \to v \quext{weakly in }\, L^p(0,T;W^{2,p}(\Omega)),
   \quext{for some }\, p>1,
\end{equation}
and $\dn v=0$ on $(0,T)\times \Gamma$ in the sense of traces.
Coupling \eqref{conv22} and \eqref{conv32} we obtain
strong convergence of $\nabla v\dd$ by interpolation. 
Indeed, we can use (for example) the Gagliardo-Nirenberg
inequality (see again \cite{Ni}) in the form
\begin{equation}\label{conv33}
  \big\| \nabla (v\dd - v) \big\|_{L^{24/19}(\Omega)}
   \le c \big\| D^2 (v\dd - v) \big\|_{L^1(\Omega)}^{1/2}
              \| v\dd - v \|_{L^{12/7}(\Omega)}^{1/2}
  + \| v\dd - v \big\|_{L^1(\Omega)}.
\end{equation}
Then, integrating in time, and using the boundedness of
$W^{2,p}$-norms resulting from \eqref{conv32}, we infer
\begin{align} \no
  & \big\| \nabla (v\dd - v) \big\|_{L^1(0,T;L^{24/19}(\Omega))}
   \le c \big\| D^2 (v\dd - v) \big\|_{L^1((0,T)\times\Omega)}^{1/2}
              \| v\dd - v \|_{L^1(0,T;L^{12/7}(\Omega))}^{1/2}
   + \| v\dd - v \big\|_{L^1((0,T)\times \Omega)}\\
 \no
  & \mbox{}~~~~~
   \le c \big( \| D^2 v\dd \|_{L^1((0,T)\times\Omega)}^{1/2}
   + \| D^2 v \|_{L^1((0,T)\times\Omega)}^{1/2} \big)
              \| v\dd - v \|_{L^1(0,T;L^{12/7}(\Omega))}^{1/2}
  + \| v\dd - v \big\|_{L^1((0,T)\times \Omega)}\\
  \label{conv34}
  & \mbox{}~~~~~
    \le c \| v\dd - v \|_{L^1(0,T;L^{12/7}(\Omega))}^{1/2}
    + \| v\dd - v \big\|_{L^1((0,T)\times \Omega)}.
\end{align}
Owing to \eqref{conv22}, we then obtain that
$\nabla v\dd$ tends to $\nabla v$, say, strongly
in $L^1((0,T)\times \Omega)$, and consequently
almost everywhere. Thus, recalling \eqref{conv23} and using
once more the generalized Lebesgue theorem, we have 
\begin{equation}\label{conv23b}
  \nabla v\dd \to \nabla v \quext{strongly in }\, L^p((0,T)\times\Omega),
   \quext{for some }\, p>2.
\end{equation}
whence, by \eqref{conv21}, \eqref{conv24} is improved up to
\begin{equation}\label{conv24b}
  \frac{u\dd}2 |\nabla v\dd|^2 \to \frac{u}2 |\nabla v|^2 
   \quext{strongly in }\, L^p((0,T)\times\Omega),
   \quext{for some }\, p>1.
\end{equation}
Consequently, the function $\Phi$ in \eqref{CH2vlim} 
is identified to its expected limit. Hence, we get
\eqref{CH2w} (holding as a relation in $L^p(\Omega)$ for some
$p>1$, hence almost everywhere) and the boundary condition 
\eqref{neumw}. Finally, we notice that the limit $\delta\searrow 0$
can be taken trivially in the initial condition \eqref{iniz-intro}. 
This concludes the proof.
\beos\label{rem:uniq}
 Theorem~\ref{teo:main} states that uniqueness
 holds for weak solutions that are classical on
 all intervals $(\tau,T)$, $\tau>0$. This is, in fact,
 the same regularity class where we are able to 
 prove existence. However, we cannot exclude 
 that uniqueness might instead fail as one considers 
 the (larger) class of {\sl all}\/ weak solutions, 
 which in particular may contain some trajectory that
 does not achieve the ``classical'' regularity for strictly positive
 times. Actually, some of the calculations given
 in the proof (in particular, those related to 
 the gradient terms) seem not be justified under 
 the sole regularity conditions proper of weak solutions 
 (which, for instance, may not be ``separated'' from 
 the singular values $\pm1$).
\eddos



\vspace{15mm}

\noindent%
{\bf First author's address:}\\[1mm]
Giulio Schimperna\\
Dipartimento di Matematica, Universit\`a degli Studi di Pavia\\
Via Ferrata, 1,~~I-27100 Pavia,~~Italy\\
E-mail:~~{\tt giusch04@unipv.it}

\vspace{4mm}

\noindent%
{\bf Second author's address:}\\[1mm]
Irena Paw\l ow\\
Systems Research Institute,\\
Polish Academy of Sciences\\
\mbox{}~~and Institute of Mathematics and Cryptology,\\
Cybernetics Faculty,\\
Military University of Technology,\\
S.~Kaliskiego 2,~~00-908 Warsaw,~~Poland\\ 
E-mail:~~{\tt Irena.Pawlow@ibspan.waw.pl}


\begin{thebibliography}{99}

 
 
\bibitem{AP96}
 H.W.~Alt and I.~Paw\l ow,
 {\sl On the entropy principle of phase transition models 
  with a conserved order parameter},
 Adv.\ Math.\ Sci.\ Appl.,
 {\bf 6} 
 (1996), 
 291--376. 

\bibitem{APNotes96}
 H.W.~Alt and I.~Paw\l ow,
 {\sl The Cahn-Hilliard equation for a polymer mixture},
  SFB 256 Universit\"{a}t Bonn, unpublished notes (1996).
  
\bibitem{Bin83}
 K.~Binder, 
 {\sl Collective diffusion, nucleation, and spinodal 
  decomposition in polymer mixtures},
 J.~Chem.\ Phys., {\bf 79} 
 (1983),
 6387--6409. 

 
 
 
\bibitem{barbu}
 V.~Barbu,
 ``Nonlinear Semigroups and Differential Equations in Banach Spaces'',
 Noord\-hoff,
 Leyden,
 1976.
 




\bibitem{BJL83}
 F.~Brochard, J.~Jouffroy, and P.~Levinson,
 {\sl Polymer-polymer diffusion in melts},
 Macromolecules, 
 {\bf 16} (1983),
 1638--1641.

 
\bibitem{Ca}
 J.W.~Cahn,
 {\sl On spinodal decomposition},
 Acta Metall., {\bf 9} (1961), 795--801.

\bibitem{CH}
 J.W.~Cahn and J.E.~Hilliard,
 {\sl Free energy of a nonuniform system. I. Interfacial free energy},
 J.\ Chem.\ Phys., {\bf 28} (1958), 258--267.
 
\bibitem{DpGG}
 R.~Dal Passo, H.~Garcke, and G.~Gr\"un,
 {\sl On a fourth-order degenerate parabolic equation: global
   entropy estimates, existence, and qualitative behavior of
   solutions},
 SIAM J.~Math.\ Anal.,
 {\bf 29} (1998), 
 321--342.
 
\bibitem{dG80}
 P.G.~de Gennes,
 {\sl Dynamics of fluctuations and spinodal decomposition in polymer
   blends},
 J.~Chem.\ Phys.,
 {\bf 72} (1980),
 4756--4763.
 
\bibitem{dG85}
 P.G.~de Gennes,
 ``Scaling Concepts in Polymer Physics'',
 Cornell Univ.\ Press, Ithaca, 1985.
 
\bibitem{DNS}
 J.~Dolbeault, B.~Nazaret, and G.~Savar\'e,
 {\sl A new class of transport distances between measures},
 Calc.\ Var.\ Partial Differential Equations,
 {\bf 34} (2009),
 193--231.
  
\bibitem{EG}
 C.M.~Elliott and H.~Garcke,
 {\sl On the Cahn-Hilliard equation with degenerate mobility},
 SIAM J.~Math.\ Anal.,
 {\bf 27} (1996), 
 404--423.
 
\bibitem{EL91}
 C.M.~Elliott and S.~Luckhaus,
 {\sl A generalized diffusion equation for phase separation of a multi-
  component mixture with interfacial free energy},
  SFB 256 Universit\"{a}t Bonn, Preprint 195 (1991).
 
\bibitem{GNC00}
 H.~Garcke and A. Novick-Cohen, 
 {\sl A singular limit for a system of 
  degenerate Cahn-Hilliard equations}, 
 Adv.\ Differential Equations,
 {\bf 5} (2000), 401--434.
  


 
%


\bibitem{Gu}
 M.~Gurtin,
 {\sl Generalized Ginzburg-Landau and Cahn-Hilliard equations
  based on a microforce balance},
 Phys.~D,
 \textbf{92} (1996),
 178--192.
 
\bibitem{KNP} 
 N.~Kenmochi, M.~Niezg\'odka, and I.~Paw\l ow, 
 {\sl Subdifferential operator approach to the Cahn-Hilliard equation with constraint},
 J.~Differential Equations,
 {\bf 117} (1995), 
 320--356.
 
 
 
 

\bibitem{MiMa90}
 V.S.~Mitlin and L.I.~Manevich,
 {\sl Kinetically stable structures in the nonlinear theory 
 of spinodal decomposition }, 
 J.~Polymer Sci.\ Part B: Polymer Physics,
 {\bf 28} (1990), 
 1--16.
  
\bibitem{MiMaE85} 
 V.S.~Mitlin, L.I.~Manevich, and I. Ya~Erukhimovich,
 {\sl Formation of kinetically stable domain structure during
  spinodal decomposition of binary polymer mixtures},
 Zh.~Eksp.\ Teor.\ Fiz.,
 {\bf 88} (1985), 495--506; 
 Sov.\ Phys.\ JETP,
 {\bf 61} (1985), 
 290--296.
   
 
\bibitem{MZ} 
 A.~Miranville and S.~Zelik, 
 {\sl Robust exponential attractors for Cahn-Hilliard type equations 
  with singular potentials},
 Math.\ Methods Appl.\ Sci.,
 {\bf 27} (2004), 
 545--582.
 
\bibitem{NB89} 
 E.B.~Nauman and N.P.~Balsara,
 {\sl Phase equilibria and the Landau-Ginzburg functional},
 Fluid Phase Equilibria,
 {\bf 45} (1989),
 229--250. 

\bibitem{NL84}
 A.E.~Nesterov and J.S.~Lipatov,
 ``Thermodynamics of Solutions and Mixtures of Polymers'',
 Naukova Dumka, Kiev, 1984 (in Russian).
 
\bibitem{Ni}
 L.~Nirenberg,
 {\sl On elliptic partial differential equations},
 Ann.\ Scuola Norm.\ Sup.\ Pisa (3),
 {\bf 13} (1959),
 115--162.


\bibitem{Nos87}
 T.~Nose,
 {\sl Kinetics of phase separation in polymer mixtures},
 Phase Transitions,
 {\bf 8} (1987),
 245--260.



\bibitem{PZ11}
 I.~Paw\l ow and W.~Zaj\c aczkowski,
 {\sl A sixth order Cahn-Hilliard type equation arising in
   oil-water-surfactant mixtures},
 Comm.\ Pure Appl.\ Anal.,
 {\bf 10} (2011), 
 1823--1847. 

\bibitem{Pin81}
 P.~Pincus, 
 {\sl Dynamics of fluctuations and spinodal decomposition 
 in polymer blends. II},
 J.~Chem.\ Phys.,
 {\bf 75} 
 (1981),
 1996--2000.

  

  
\bibitem{SP11}
 G.~Schimperna and I.~Paw\l ow,
 {\sl On a Cahn-Hilliard model with nonlinear diffusion},
 preprint arXiv:1106.1581 (2011),
 submitted.  


\bibitem{Si}
 J.~Simon,
 {\sl Compact sets in the space {$L^p(0,T;B)$}},
 Ann.\ Mat.\ Pura Appl.\ (4),
 {\bf 146} (1987),
 65--96.

 

\bibitem{Wit98}
T.P.~Witelski,
 {\sl Equilibrium interface solution of a degenerate singular 
  Cahn-Hilliard equation},
 Appl.\ Math.\ Lett.,
 {\bf 11} (1998), 
 127--133.

\bibitem{ZZE06}
D.~Zhou, P.~Zhang, and W.~E, 
{\sl Modified models of polymer phase separation},
 Phys.\ Rev.~E,
 {\bf 73} (2006),
 061801.


\end{thebibliography}
\end{document}